\documentclass[11pt]{amsart}
\usepackage{amssymb}

\newtheorem{prop}{Proposition}
\newtheorem{lem}{Lemma}
\newtheorem{thm}{Theorem}
\newtheorem{cor}{Corollary}

\newcommand{\codim}{\mbox{codim}}
\newcommand{\proj}{\mathbf P}

\newcommand{\rarr}{\rightarrow}
\newcommand{\oh}{{\mathcal{O}}}
\newcommand{\f}{{\mathcal{F}}}
\newcommand{\com}{\mathbb{C}}
\newcommand{\Z}{\mathbb{Z}}

\newcommand{\Pone}{{\proj^1}}
\newcommand{\F}{\mathbf F}
\newcommand{\flag}{{\F(n)}}
\newcommand{\pflag}{\F^{N}}
\newcommand{\hquot}{\mathcal{HQ}_\bd(\F(n))}
\newcommand{\hquote}{\mathcal{HQ}_{\bd-\be}}
\newcommand{\hquotei}{\mathcal{HQ}_{\bdei}}
\newcommand{\hquotem}{\mathcal{HQ}_{\bdem}}

\newcommand{\D}{\mathcal{D}}
\newcommand{\e}{\tilde{e}}
\newcommand{\E}{\tilde{E}}
\newcommand{\bel}{{\mathbf{e}_l}}
\newcommand{\Rl}{\mathcal{HQ}_\bel}
\newcommand{\oO}{\overline{\Omega}}

\newcommand{\bc}{\mathbf{c}}
\newcommand{\bd}{\mathbf{d}}
\newcommand{\be}{\mathbf{e}}

\newcommand{\bdem}{\bd-\bem}
\newcommand{\bem}{{\mathbf{e}_m}}
\newcommand{\bde}{{\bd-\be}}
\newcommand{\bei}{{\mathbf{e}_i}}
\newcommand{\bdei}{{\bd-\bei}}
\newcommand{\bq}{\mathbf{q}}

\newcommand{\BQ}{\mathbf{Q[1]}}
\newcommand{\qsum}{\sum_{\bd,w}\bq^\bd}
\newcommand{\qqsum}{\sum_{\bd,w'}\bq^\bd}
\newcommand{\sw}{\sigma_{\dualw}}
\newcommand{\sww}{\sigma_{\dualww}}
\newcommand{\dualw}{{w^{\vee}}}
\newcommand{\dualww}{{{w'}^{\vee}}}
\newcommand{\fS}{\mathfrak{S}}
\newcommand{\U}{\mathcal{U}}

\newcommand{\R}{\mathcal{HQ}_\bd}
\newcommand{\Rde}{\mathcal{HQ}_\bde}

\newcommand{\dd}{{d_1,\ldots,d_{n-1}}}
\newcommand{\hrarr}{\hookrightarrow}
\newcommand{\A}{\mathcal{A}}

\newcommand{\Q}{\mathcal{Q}}

\newcommand{\rrarr}{\twoheadrightarrow}

\newcommand{\eqq}{\stackrel{\sim}{=}}

\newcommand{\bpf}{\noindent {\em Proof.} }
\newcommand{\epf}{\qed \vspace{+10pt}}

\newcommand{\mor}{\mathrm{Mor}}

\begin{document}
\title{Quantum Cohomology of Flag Manifolds}
\author{Linda Chen}

\maketitle

\section{Introduction}

In this paper, we study the (small)
 quantum cohomology ring of the partial flag manifold.
We give proofs of the presentation of the ring and
of the quantum Giambelli formula for Schubert varieties.  These are 
known results, but our proofs are more natural and direct 
than the previous ones.  

One of our goals is to give evidence of a
relationship between {\it universal Schubert polynomials}, which give the
answer to a degeneracy locus problem, and
{\it quantum Schubert polynomials}, which appear in quantum cohomology.
It has been known that the universal Schubert polynomials
specialize to both the ordinary and the quantum Schubert polynomials, but
previous reasons for this have been purely algebraic \cite{fu}.

The quantum cohomology ring of a projective manifold $X$
is a deformation of the ordinary cohomology ring of $X$.
The classical Schubert calculus, consisting of Giambelli and Pieri-type 
formulas which give the multiplicative 
structure of the cohomology ring of the flag manifold, 
has been used as a tool to solve enumerative problems.  Similarly,
the entries in the quantum
multiplication table count rational curves  on a flag manifold of a given
multidegree which meet three general Schubert varieties. 
These numbers can be interpreted 
as intersection numbers on appropriate moduli spaces of holomorphic maps from
the projective line $\Pone$ to the flag manifold.

In order to understand these intersections, various compactifications
of the moduli space of maps have been studied, 
for example the stable maps of Kontsevich. However, in the case of partial flag
manifolds, including Grassmannians and complete flag manifolds, there
are smooth compactifications called {\it hyperquot schemes}, which
 generalize Grothendieck's Quot scheme \cite{g}.  They have been studied by 
Ciocan-Fontanine \cite{cf1} \cite{cf2}, Laumon \cite{lau},
Kim \cite{k}, and in \cite{c}. 
Most of what is known about
the quantum cohomology of flag manifolds rely heavily on computations in
the cohomology of hyperquot schemes. We obtain our results through a further
study of the intersection theory of  hyperquot schemes.

For the sake of notation, we first state and prove our results 
for the case of the complete flag manifold.  Most of the statements hold
verbatim for the general case of partial flag manifolds, and 
many of the proofs need only slight modifications.  We give
ingredients to extend the arguments in the final section of this paper.  
The exceptions to this are found in
sections \ref{other} and \ref{proof}, whose constructions and results 
apply only to complete
flag manifolds, and whose methods do not generalize.  
In section \ref{qcpf}, we give an alternate approach completely
bypasses this special argument.  
We include sections \ref{other} and \ref{proof} because they may be of
outside interest as we introduce and study a new set of degeneracy
loci on the hyperquot scheme.

Two components of a classical {\it Schubert calculus} are a presentation
of the cohomology ring and a Giambelli formula, which writes Schubert
classes in terms of the generators of the ring.
The classical cohomology ring of the flag manifold $\F$
has a presentation $\Z[\sigma_{i}^j]/I$, where 
the $\sigma_{i}^j$'s are determined
by Chern classes of certain tautological vector bundles on $\F$,
and the relations are given by the ideal $I$.
By the general result of Siebert and Tian in \cite{st}, 
the quantum cohomology ring of $\F$ has presentation
$\Z[\sigma_{i}^j,q_1,\ldots,q_m]/I_q$, where 
the variables $q_i$ are deformation parameters and $I_q$ is
a deformation of the ideal $I$.

The classical Giambelli formula gives Schubert classes as polynomials
in the variables $\sigma_i^j$ called {\it Schubert polynomials}.  
{\it Special Schubert classes} are those corresponding to Chern classes
of the tautological vector bundles.
The quantum Giambelli formula is a deformation of the classical formula,
 giving the
Schubert classes as polynomials in the variables $\sigma_i^j$ and  $q_k$.

For complete flag manifolds, 
generators of the ideal $I_q$ were conjectured by
Givental and Kim in \cite{gk}.  This and
the special case of the quantum Giambelli formula were simultaneously
proved by Ciocan-Fontanine in \cite{cf1}.  Using
these results and combinatorial methods, Fomin, Gelfand, and Postnikov
constructed {\it quantum Schubert polynomials}, and proved the general case
of the quantum Giambelli formula \cite{fgp}.  An independent proof of the
presentation of the ring was given by Kim in \cite{k}.  The methods of
\cite{cf1} and \cite{fgp} were adapted to 
prove the results for partial flag manifolds \cite{cf2}.

We provide a simplified argument for both of these results, 
which simultaneously proves the presentation and the quantum Giambelli formula.
One of the key points in the proof is to
use a certain degeneracy locus formula for flags of bundles, stated 
 in terms of universal Schubert polynomials, and proved by Fulton
 \cite{fu}.  Our methods of proof are
similar to those
used by Bertram  and Ciocan-Fontanine, including
a liberal use of certain maps constructed by Ciocan-Fontanine to 
understand the boundary of the hyperquot scheme \cite{b}\cite{cf1}\cite{cf2}.

In sections \ref{hyperquot} and \ref{schubert}, we review the main 
constructions of \cite{cf1} which are used in the paper.   
In section \ref{univsb}, we review the degeneracy locus formula of Fulton
as it applies to our situation. In section \ref{results}, we prove
the results of the paper via a Main Proposition, which is proved in section
\ref{mainprop}, modulo a lemma stated as Proposition \ref{q1mult}. 
Its general statement for partial flag manifolds
is stated and proved in section \ref{qcpf}. 
In the complete flag manifold case, an alternate proof of the lemma
relies on the
new ideas and constructions found in sections \ref{other} and \ref{proof}.

{\em Acknowledgements.} This paper gives a generalization of part of my
Ph.D. thesis. The methods originate in the work of  Aaron Bertram and
Ionu\c{t} Ciocan-Fontanine  on the quantum cohomology of
Grassmannians and flag manifolds.  Their ideas and constructions have
inspired the approach. I thank William Fulton 
for his encouragement and helpful comments, and for introducing me to the
problem.  This research was partially supported by the Clay Mathematics
Institute.

\section{Classical cohomology of the flag manifold}
\label{flag}

Let $\flag$ denote the complete flag manifold of $\com^n$,
which parametrizes flags of subspaces $V_\bullet$:    
$$\{0\}=V_0\subset V_1 \subset V_2 \subset ... \subset V_{n-1} 
\subset V=\com^n$$
with dim $V_i=i$.  Write $V_X=V\otimes \oh_X$ for any scheme $X$.

There is a universal sequence of vector bundles on $\flag$:
$$V_{\flag}\rrarr Q_{n-1}\rrarr\cdots\rrarr Q_{1}$$
with $\mbox{rank }Q_i=i$, and each $Q_i\rarr Q_{i-1}$ a surjection. 
Any flag of successive vector bundle quotients of $V_X$:
$$V_X\rrarr F_{n-1}\rrarr\cdots\rrarr F_{1}$$
with $\mbox{rank }F_i=i$ gives a map $f:X\rarr\flag$ with
$F_i=f^*Q_i$ where the $Q_i$'s are the 
tautological quotient bundles on $\flag$.

Fix a flag $V_\bullet: V_1\subset \cdots V_{n-1}\subset V_n$.  
For any $w$ in the symmetric group $S_n$, define
$r_w(q,p)$ to be the number of $i\leq q$ such that $w(i)\leq p$.
The corresponding {\it Schubert variety} is given by:
$$\Omega_w(V_\bullet)=\{U_\bullet\in\flag: \mbox{rank}_{U_\bullet}(V_p\otimes 
\oh_{\flag}\rarr Q_q)\leq r_w(q,p)\mbox{ for }1\leq p,q \leq n\}.$$
This is a codimension $l(w)$ subvariety in $\flag$, where $l(w)$ is the length
of the permutation $w$.  Its class is independent of choice of flag
$V_\bullet$, and we denote this class by $[\Omega_w]$.  

Let $w_0\in S_n$ be the permutation of longest length,  with $w_0(i)=n-i+1$
for $1\leq i\leq n$.  For $w\in S_n$, write $\dualw=w_0w$.
We have the following classical results.
\begin{thm} The classes $[\Omega_w]$ form an additive basis for
$H^*(\F)$.  Furthermore, for $w\in S_n$, the Schubert classes
$[\Omega_w]$ and $[\Omega_{\dualw}]$ are Poincar\'{e} dual.
\end{thm}

\begin{thm} Let $x_i=c_1(Q_{i}\rarr Q_{i-1})$ 
for $1\leq i\leq n$.
Then $$H^*(\flag,\Z)\eqq \Z[x_1,\ldots,x_n]/(e_1(n),\ldots,e_n(n))$$
where $e_i(n)$ is the $i$th symmetric polynomial in $x_1,\ldots,x_n$.
\end{thm}

The Giambelli problem is to express $[\Omega_w]$ in terms of this presentation.
To do this, we give the definition of {\it Schubert polynomials} as
given by Lascoux and Schutzenberger \cite{ls}.  
For $1\leq i \leq n-1$, let $\partial_i$
act on $\Z[x_1,\ldots,x_n]$ by
$$\partial_i P= \frac{P(x_1,\ldots,x_n)-P(x_1,\ldots,x_{i-1},x_{i+1},x_i,x_{i+2},\ldots,x_n)}{x_i-x_{i+1}}.$$

Let $s_i$ be the transposition $(i,i+1)$.  For
$w\in S_n$, write $w=w_0\circ s_{i_1}\circ\cdots\circ s_{i_k}$, where
$k=\binom{n}{2}-l(w)$.  Then the
Schubert polynomial associated to $w$ is defined by
$$\fS_w(x)= \partial_{i_k}\circ\cdots\circ\partial_{i_1}
(x_1^{n-1}x_2^{n-2}\cdots x_{n-1}).$$

The solution to the Giambelli problem was given by Bernstein,
Gelfand, and Gelfand \cite{bgg} and Demazure \cite{d}, cf. \cite{mac}.
\begin{thm} $[\Omega_w] = \fS_w(x)$ in $H^*(\flag,\Z)$.
\end{thm}

\section{Quantum multiplication map}
\label{qmult}

Additively, we can view
the small quantum cohomology ring as $$QH^*(\flag)=H^*(\flag)\otimes_\Z
\Z[q_1,\ldots,q_{n-1}].$$ Then as $\Z[q]-$modules, the 
ordinary cohomology has a canonical injection into
the quantum ring given by $[\Omega_w]\rarr[\Omega_w]\otimes 1$ .  
Indeed, $QH^*(\flag)$ is a deformation of $H^*(\flag)$
so that the ordinary ring is recovered by setting $q_i=0$.
 
For $w\in S_n$, 
let $\sigma_w:=[\Omega_w]\otimes 1$ denote the Schubert class in the small
quantum cohomology ring, so that substituting
$q_i=0$ into $\sigma_w$ gives $[\Omega_w]=\fS_w(x)$, the 
Schubert class defined in section \ref{flag}. 
(This notation 
differs from that used by
Bertram and Ciocan-Fontantine \cite{b} \cite{cf1} \cite{cf2}.  We use
the variable $\mu$ instead of $\sigma$ to distinguish between the
class $\sigma_w$ in quantum cohomology of the flag manifold and the class 
$\mu_w$ in the cohomology of the hyperquot scheme, which is defined in
section \ref{schubert}.)

Let $s_i\in S_n$ be the transposition $(i,i+1)$.
We say that a map $f:\Pone\rarr \flag$ has multidegree $\bd= (d_1,...,d_{n-1})$
when $f_*[\Pone] = \sum d_i[\Omega_{w_0s_i}]$, with each $d_i$ a nonnegative 
positive integer.  Recall that
$[\Omega_{w_0s_i}]$ is dual to $[\Omega_{s_i}]$.

The  {\it Gromov-Witten number}
\begin{equation}
\label{GWnumber}
 \langle \Omega_{w_1},\ldots,\Omega_{w_N}\rangle_\bd \in\Z
\end{equation}
is defined as follows.  For any $t_1,\ldots,t_N\in\Pone$ in general position, 
and $\Omega_{w_1},\ldots,\Omega_{w_N}\subset \flag$
 general translates (obtained by choosing general flags), 
it is 
the number of holomorphic maps $f:\Pone\rarr\flag$ of multidegree 
$\bd$
satisfying $f(t_i)\in\Omega_{w_i}$ for $1\leq i\leq N$ if this number is
finite, and zero otherwise.

We write $\bq^\bd = q_1^{d_1}q_2^{d_2}\cdots q_{n-1}^{d_{n-1}}$.
Define the quantum multiplication map:
\begin{equation}
\label{qmultdef}
\sigma_{w_1}* \sigma_{w_2} =
\qsum \langle \Omega_{w_1},\Omega_{w_2},\Omega_w \rangle_\bd \sw.
\end{equation}
This product 
gives $QH^*(\flag)$ the structure of a commutative, associative 
$\Z[q]-$algebra.  This follows from general associativity results on
the big quantum ring which can be specialized 
to the small quantum ring, see \cite{fp}.

With this defined multiplication, we have:
\begin{prop}
\label{multmap}
$$\sigma_{w_1}*\ldots * \sigma_{w_N} =
\qsum \langle \Omega_{w_1},\ldots,\Omega_{w_N},\Omega_w \rangle_\bd \sw.$$
\end{prop}
This was first proved for Grassmannians \cite{b} by using the Quot scheme.
It has been proved for Grassmannians \cite{fp}  and for partial flag manifolds
 \cite{cf2} by realizing  the Gromov-Witten numbers defined in (\ref{GWnumber})
as intersection numbers on $\overline{M}_{0,N+1}(\flag,\bd)$. 
The validity of this statement of the quantum product is also
proved via the hyperquot scheme and 
results of section \ref{schubert} in \cite{cf1}, using the methods of Bertram.

\section{The hyperquot scheme}
\label{hyperquot}

In this section, we summarize the constructions and results that we need for
our proofs.  In particular, we describe the hyperquot scheme, some degeneracy
loci, and a description of the boundary.  The proofs of the results of
section \ref{boundary} are
found in  detail in \cite{cf1}. 

\subsection{Properties of $\R$}
Consider the following functor $\f_\bd$ from the category of schemes to the 
category of sets.  For a scheme $T$, 
$\f_\bd (T)$ is defined to be the set of equivalence classes of 
flagged quotient sheaves
$$V^*_{\Pone \times T}\rrarr \Q_{n-1} \rrarr \cdots \rrarr \Q_{1}$$
with each $\Q_i$ flat over $T$ with Hilbert polynomial
$\chi(\Pone_t,(\Q_i)_t(m)) = (m+1)i+d_{n-i}$ on the fibers
of $\pi_T: \Pone \times T \rightarrow T$, so that $\Q_i$ 
is of rank $i$ and relative degree $d_{n-i}$ over $T$,
i.e. that $(\Q_i)_t$ is of degree $d_{n-i}$ for every $t\in T$.  
Two such flags 
$V^*_{\Pone \times T}\rrarr \Q_{n-1}^1 \rrarr \cdots \rrarr \Q_1^{1}$ 
and $V^*_{\Pone \times T}\rrarr \Q_{n-1}^2 \rrarr \cdots \rrarr \Q_1^{2}$
are in the same equivalence class when there exist maps $\Q^1_i\rarr Q^2_i$
so that all squares commute.

The functor $\f_\bd$ is represented by the projective
scheme $\R=\hquot$ \cite{cf1} \cite{k}.  Its construction as  the fine moduli space of flat families of flagged 
quotient sheaves over $\Pone$
has been described by Ciocan-Fontanine following the
ideas of Grothendieck and Mumford \cite {cf1} \cite{g} \cite{m}. 
It has also been 
described in a different way by Kim \cite{k}, 
as a closed subscheme of a product of 
Quot schemes, and independently by Laumon \cite{lau}. We have

\begin{thm}
\label{Rpfproj}
$\hquot$ is an irreducible, rational, nonsingular, projective variety
of dimension  $\binom{n}{2}+2\sum d_i$.
 \end{thm}

Thus, associated to  $\hquot$
is a universal sequence of sheaves on $\Pone \times \R$ of
successive quotients of sheaves, each of which is flat over $\R$:
$$ V^*_{\Pone\times \R} 
\twoheadrightarrow B_{n-1} \twoheadrightarrow \cdots \twoheadrightarrow B_{1}.$$
In general, the sheaf $B_i$ is not locally free.  
Consider the sheaves
 $A_i:=\ker(V^*_{\Pone\times \R}\rarr B_{n-i})$.  
Each $A_i$ is flat over $\R$, and it
is an easy consequence of flatness and the fact that $\Pone$ is a 
nonsingular curve that each $A_i$ is locally free.
Thus, we have the following universal sequence on  $\Pone \times \R$:
\begin{equation}
\label{univ} 
A_1 \hrarr A_2 \hrarr \cdots \hrarr A_{n-1} \hrarr V^*_{\Pone\times \R} 
\twoheadrightarrow B_{n-1} \twoheadrightarrow \cdots \twoheadrightarrow B_{1}.
\end{equation}
with $A_i$ and $B_i$ of rank $i$.  Denote the
inclusion maps by $\gamma_i: A_i \hrarr A_{i+1}$ and the surjections
by $\pi_i:B_{i+1} \rarr B_{i}$ for each $1\leq i\leq n-1$.
We set  $A_{n}=B_n= V^*_{\Pone\times \R}$ and
$A_0=B_{0}=0$.  
The map  $\gamma_i: A_i\hrarr A_{i+1}$ is an inclusion of sheaves, 
not an inclusion of bundles.

Consider $\mor_\bd(\Pone, \flag)$, the space of morphisms
from $\Pone$ to $\flag$
 of multidegree $\bd = (d_1,...,d_{n-1})$.
By the universal property of $\flag$, a morphism $f \in\mor_\bd(\Pone, \flag)$
corresponds to 
successive quotient bundles
$$V_\Pone \rarr Q_{n-1} \rarr \cdots \rarr Q_1$$
with $Q_i$ of rank $i$ and degree $d_{n-i}$. 
Equivalently, by taking kernels, $f$ corresponds to successive
subbundles $S_1\hrarr \cdots \hrarr S_{n-1}\hrarr V_\Pone$. 
Dualizing, we see that
$\mor_\bd(\Pone, \flag)$ parametrizes  successive quotient bundles 
of $V^*_\Pone$ of rank $i$ and degree $d_{n-i}$.  
In this way, the hyperquot scheme
$\R$ is a compactification of $\mor_\bd$.  

We define the hyperquot
scheme in terms of the dual trivial vector bundle
 $V^*_\Pone$ instead of $V_\Pone$ to
ensure that the universal subsheaves $A_i \hrarr V^*_\Pone$ are locally
free.  This allows us to apply the degeneracy locus formula of 
section \ref{deg} more readily.

\subsection{Description and construction of $\U_\be$}
\label{constructU}
In the next two sections, we review the construction of 
certain schemes which map to various hyperquot schemes \cite{cf1}.  
These morphisms are  used 
to understand a recursive structure of the boundary of $\R$, which 
 is used to understand intersections of degeneracy loci 
on the hyperquot scheme, to be defined in sections \ref{schubert} and
\ref{other}.

Let $\be=(e_1,\ldots,e_{n-1})$ be a sequence of nonnegative integers
satisfying:
\begin{enumerate}
\item[(1)] $e_i\leq \min(i,d_i)$ for $1\leq i\leq n-1$,
\item[(2)] $e_i-e_{i-1}\leq 1$ for $2\leq i\leq n-1$.
\item[(3)] $\sum e_i \geq 1$
\end{enumerate}

We prove a lemma to be used in section \ref{other}.
\begin{lem}
\label{es}
For $\be$ satisfying $(1),(2),$ and $(3)$, and setting $e_0$=0,
\begin{enumerate}
\item $\sum e_i \leq \sum e_i(1+e_i-e_{i-1})$
\item $\sum e_i(1+e_i-e_{i-1})\geq 2$ with equality if and only if $\sum e_i =1$.
\end{enumerate}
\end{lem}
\bpf This follows from the observation that
 $$\sum e_i(1+e_i-e_{i-1}) = 
\sum e_i +\frac{1}{2}[e_1^2+(e_2-e_1)^2+\cdots 
(e_{n-1}-e_{n-2})^2+e^2_{n-1}].$$
\epf

For each such multiindex $\be$, we consider the scheme $\U_\be$ as
in \cite{cf1} as follows.
On $\Pone\times\Rde$, there is the universal sequence
$$A_1^\bde\hrarr \cdots \hrarr A_{n-1}^\bde\hrarr V^*_{\Pone\times\Rde}.$$

For $1\leq i\leq n-1$, let
$X_i$ be the Grassmann bundle of $e_i$-dimensional quotients of 
$A_i^\bde$, and let $X_\be$ be the fiber product of these $X_i$'s,
with projection map $\pi:X_\be \rarr \Pone\times\Rde$.  Let $K_i$ denote
the tautological subbundle of the pullback of $A_i^\bde$
over $X_i$, and $Q_i$ the corresponding quotient bundle, so that for 
$1\leq i\leq n-1$ there is the exact sequence
$$0\rarr K_i\rarr A_i^\bde \rarr Q_i$$
with $K_i$ of rank
$i-e_i$ and $Q_i$ of rank $e_i$.  Let $K_i$ and $Q_i$ also denote
the pullbacks of these bundles to $X_\be$ via the natural projections.  

Define $\U_\be$ to be the locally closed subscheme of $X_\be$ given
 by the conditions
\begin{enumerate}
\item[(1)] The composite map $K_i\rarr A_i\rarr A_{i+1}\rarr Q_{i+1}$ vanishes 
for $1\leq i\leq n-2,$
\item[(2)] $\mbox{rank }(K_i\rarr V^*_{X_\be})=i-e_i$ for $1\leq i\leq n-1.$
\end{enumerate}

We review an explicit construction of $\U_\be$, as well as the construction 
of the morphism $h_\be:\U_\be\rarr\R$, 
which we  use in later sections.
Let $U\subset X_\be$ be the open subscheme given conditions (2).  
We construct $\U_\be$ inductively as a sequence of Grassmann bundles. 
\begin{enumerate}
\item Let $\rho_1:\U_1:= G^{e_1}(\pi^*A_1^{\bd-\be})\rarr U$ 
be the Grassmann bundle of $e_1$-dimensional
quotients, with universal 
subbundle $K_1$.

\item Consider the Grassmann bundle 
$\rho_i: \U_i := G^{e_i}(\rho^*\pi^*A_i^{\bd-\be}/K_{i-1})\rarr \U_{i-1}$, 
with universal subbundle $S_i$, where $\rho$ denotes the composition
$\rho_{i-1}\circ\cdots\circ\rho_1$.
Let $K_i$ be the natural extension of $\rho_i^*K_{i-1}$ by $S_i$.
\end{enumerate}
Then $\U_\be=\U_{n-1}$.

\subsection{Construction of $h_\be$}
Let $\pi:\U_\be\rarr\Pone\times\Rde$ be the projection map.
Consider the map 
$$\psi:= 1\times \pi:\Pone\times\U_\be \rarr \Pone\times\Pone\times\Rde.$$
Let $\Delta\subset \Pone\times\Pone$ be the diagonal, and 
$\tilde{\Delta}\subset\Pone\times\U_\be$ be $\psi^{-1}(\Delta\times\Rde)$.
Denote by $p:\Pone\times\U_\be\rarr\U_\be$ the second projection.
For $1\leq i\leq n-1$, define $\widetilde{A}_i^\be$ to be the kernel of the
map $$p^*\pi^*A_i^\bde \rarr p^*Q_i|_{\tilde{\Delta}}.$$

This gives a sequence of sheaf injections
$$\widetilde{A}_1^\be\hrarr\cdots\hrarr\widetilde{A}_{n-1}^\be$$
with each $\widetilde{A}_i^\be$ flat over $\U_\be$,
locally free of rank $i$, and of
relative degree $-d_i$ on $\Pone\times\U_\be$.  Since $\R$ represents the hyperquot functor,
this defines a morphism $h_\be:\U_\be\rarr\R$, which satisies 
$(1\times h_\be)^*A_i^\bd = \widetilde{A}_i^\be$ for every $1\leq i\leq n-1$.

\subsection{The boundary}
\label{boundary}

\begin{thm}\cite{cf1} 
\label {Ue} 
Let the multiindex $\be=(e_1,\ldots,e_{n-1})$, $\U_\be$, and $h_\be$
 be as above.  Then

\begin{enumerate}  \item $\U_\be$ is smooth, irreducible, and of dimension
$$ \binom{n}{2}+2|\bd| +1 -\sum_{i=1}^{n-1}e_i(1+e_i-e_{i-1}).$$  The projection map $\pi:\U_\be\rarr\Pone\times\Rde$
is smooth and proper, with irreducible fibers.

\item If $\mbox{rank}_{(t,x)}B_i^\bd=n-i+e_i$ for $1\leq i\leq n-1$
at $(t,x)\in\Pone\times\R$, then $x\in h_\be(\U_\be)$.

\end{enumerate}
\end{thm}
 
Part 2 of the theorem implies that the boundary of the hyperquot scheme
$\R\smallsetminus\mor_\bd$ is covered by the images of the $\U_\be$ under
the morphisms $h_\be$.  

For any $t\in\Pone$ and any multiindex $\be$ as above, 
define $\U_\be(t) = \pi^{-1}(t\times\hquote)$.

\section{Schubert varieties on the hyperquot scheme}
\label{schubert}
We describe subschemes of the hyperquot scheme which are determined by the
degeneracy conditions which describe Schubert varieties, 
whose associated intersection numbers
are equal to certain Gromov-Witten invariants.

Consider the evaluation map $$ev:\Pone\times\mor_\bd\rarr \flag$$
given by $ev(t,f)=f(t)$ and use this to pull back Schubert varieties of
$\flag$ to $\mor_\bd$ in the following manner. Define for any $w\in S_n$,
$$\Omega_w(t)=ev^{-1}(\Omega_w)\cap (\{t\}\times\mor_\bd).$$

We wish to extend $\Omega_w(t)$ over the boundary
to the entire hyperquot scheme.  This can be achieved as a
certain degeneracy locus.  We fix a flag $V_\bullet$.
On $\Pone\times\hquot$,  we have the dualized
universal sequence of subsheaves
\begin{equation*}
V_{\Pone\times\R}\rarr A_{n-1}^* \rarr\cdots\rarr A_1^*.
\end{equation*}
Define $\overline{\Omega}_w$ to be the locus where $$\mbox{rank }
(V_p\otimes \oh_{\Pone\times\R}\rarr A_q^*)\leq r_w(q,p),
1\leq p,q\leq n\},$$
with the natural scheme structure given by 
vanishing of determinants,
and define $\overline{\Omega}_w(t)$ to be its  restriction to 
$\{t\}\times \R$, viewed as a subscheme of $\R$ via the identification
$t\times\R \eqq \R$.

\begin{lem}
\label{hO} We can write
$$h_\be^{-1}(\overline{\Omega}_w(t))= \pi^{-1}(\Pone\times \overline{\Omega}_w(t)) 
\cup \widetilde{\Omega}_w(t)$$
with $\widetilde{{\Omega}}_w(t)$ the degeneracy locus inside 
$\pi^{-1}(t\times \Rde)$
given by $$\mbox{rank }(V_p\rarr K_q^*)\leq r_w(q,p).$$  This equality is
scheme-theoretic away from the intersection.
\end{lem}

\begin{lem}
\label{codimO}
For $w\in S_n$
$$\mbox{codim}_{\U_\be(t)} \widetilde{\Omega}_w(t) = l(\widetilde{w}^\be)$$
for a permutation $\widetilde{w}^\be$ as described in
Construction 3.5 of \cite{cf1},
which satisfies $l(w)-l(\widetilde{w}^\be)\leq \sum e_i$.
\end{lem}

The argument of Ciocan-Fontanine in \cite{cf1}, which uses
the constructions in section  \ref{hyperquot} and
 Lemmas \ref{hO} and \ref{codimO}, 
gives the following general position result and two corollaries:

\begin{thm}
\label{gpschubert}
\begin{enumerate} \item For any subvariety $Y$ in $\R$, 
$w\in S_n$, a general translate $\Omega_w\subset \flag$, and any $t\in\Pone$,
$Y\cap \Omega_w(t)$ is either empty or pure codimension $l(w)$ in $Y$.

\item If $t_1,\ldots,t_N$ are distinct points in $\Pone$, then for general translates
of $\Omega_{w_i}$, the intersection $\bigcap_{i=1}^N \overline{\Omega}_{w_i}(t_i)$
is either empty of pure codimension $\sum_{i=1}^N l(w_i)$ in $\R$ and
is the Zariski closure of $\bigcap_{i=1}^N {\Omega_{w_i}}(t_i)$
\end{enumerate}
\end{thm}

\begin{cor} The class of $\overline{\Omega}_w(t)$ in $H^{2l(w)}(\R)$
is independent of $t\in\Pone$ and flag $V_\bullet$.  Denote this class by
$\mu_w(\bd)$.
\end{cor}
The multidegree $\bd$ is often understood, and in these cases we write
$\mu_w=\mu_w(\bd)$.  

\begin{cor}  \label{Omegaprod}
If $\sum_{i=1}^N l(w_i)= \dim\R$, and $t_1,\ldots,t_N$
are distinct points in $\Pone$, then  $\bigcap_{i=1}^N \Omega_{w_i}(t_i)=
\bigcap_{i=1}^N \overline{\Omega}_{w_i}(t_i)$, and the number of points in this
intersection is the degree of
the product $\mu_{w_1}\cdot\mu_{w_2}\cdots\mu_{w_N}$ in
$H^*(\R,\com)$.
\end{cor}

Recall our definition of the number
$ \langle \Omega_{w_1},\ldots,\Omega_{w_N}\rangle_\bd $
in section \ref{qmult}.  For $t_1,\ldots,t_N\in \Pone$ in general position,
and $\Omega_*$ general translates, this number is equal to 
the number of points of $\bigcap_{i=1}^N \Omega_{w_i}(t_i)$.  Corollary 
\ref{Omegaprod} shows that we have the following equality
\begin{equation}
\label{GWquot}
\langle \Omega_{w_1},\ldots,\Omega_{w_N}\rangle_\bd
=(\mu_{w_1}\cdot\mu_{w_2}\cdots\mu_{w_N})_\bd
\end{equation}
for every multidegree $\bd$.  The number on the left is the Gromov-Witten
number defined in (\ref{GWnumber}) and the number on the right is
the degree of the intersection product 
$\mu_{w_1}\cdot\mu_{w_2}\cdots\mu_{w_N}$ 
in the cohomology of the hyperquot scheme $\R$.

\section{Universal Schubert polynomials}
\label{univsb}
Universal Schubert polynomials, introduced in \cite{fu}, specialize to all 
known types of Schubert polynomials, including the classical version
defined in section \ref{flag} as well as the quantum Schubert polynomials
defined in \cite{fgp}.
They appear as the answer to a certain degeneracy locus problem
which we use for the proofs of our results.

\subsection{Definition of $\fS_w(c)$ and $\fS_w(g)$}
\label{usbdef}
We give two equivalent formulations of universal Schubert polynomials.
For any $w\in S_{n+1}$, the corresponding classical Schubert polynomial 
can be written
\begin{equation}
\label{spsum}
\fS_w(x)=\sum a_{k_1,\ldots,k_n}e_{k_1}(1)\cdots e_{k_n}(n)
\end{equation}
where the sum ranges over all sequences $(k_1,\ldots,k_n)$ with $k_p\leq p$
and $\sum k_i = l(w)$, and $e_k(l)$ is the $k$th elementary symmetric
polynomial in $x_1,\ldots,x_l$.  For each $w$,
the coefficients $a_{k_1,\ldots,k_n}$ are uniquely determined integers.

First, for $1\leq k\leq l \leq n$, consider independent variables $c_k(l)$.
We set $c_0(l)=1$ and $c_k(l)=0$ when $k<0$ or $k>l$.
With this notation,  the universal Schubert polynomial is defined by
$$\fS_w(c)=\sum a_{k_1,\ldots,k_n}c_{k_1}(1)\cdots c_{k_n}(n).$$

A second description uses independent variables $g_i[j]$,
where $i\geq 1, j\geq 0$, and $i+j\leq n+1$, and each $g_i[j]$ is
of degree $j+1$.  
Consider the Dynkin diagram for $(A_{n})$, mark the vertices
$x_1,\ldots,x_n$, the edges  $g_1[1],\ldots,g_{n-1}[1]$, and
the paths covering the $j+1$ consecutive vertices
$x_i,\ldots,x_{i+j}$ by $g_i[j]$.  In this notation,
we have $g_i[0]=x_i$.  
Denote by $E_k(l)(g)$ the sum of all monomials in the paths $g_i[j]$
covering exactly $k$ of the vertices $x_1,\ldots,x_l$, with
no vertex covered more than once.  If the variables $g_*[*]$ are understood,
we may sometimes write $E_k(l)$ for $E_k(l)(g)$.
Equivalently, $E_k(l)$ can defined inductively as
\begin{equation}
\label{recursion}
E_k(l)=E_k(l-1)+ \sum_{j=0}^k E_{k-j-1}(l-j-1)g_{l-j}[j].
\end{equation}
Let $G_l$ be the $l\times l$ matrix with entries $g_i[j-i]$ in the $(i,j)$ 
position for $1\leq i\leq j \leq l$, $-1$ in the $(i+1,i)$ position, 
and $0$ elsewhere. Then (\ref{recursion}) is also equivalent to  
$E_k(l)$ being the coefficient of $\lambda^k$ in the determinant
of $G_l+\lambda I$. The universal Schubert polynomial $\fS_w(g)$
is obtained by performing the substitution $c_k(l)=E_k(l)$ into $\fS_w(c)$.

Via this substution, the polynomial rings generated by $c_*(*)$ and $g_*[*]$ 
are the same, so that
the two formulations are equivalent.  As we have seen, each $c_k(l)$ can
be written as a polynomial in the $g$'s, but each $g_i[j]$ can be
written in terms of the $c_k(l)$'s as well.  In particular, 
from (\ref{recursion}), we have
\begin{equation}
\label{recursion2}
c_k(l) = c_k(l-1)+\sum_{j=0}^{k-1} c_{k-j-1}(l-j-1)g_{l-j}[j]+
g_{l-k}[k].
\end{equation}
The $g_i[j]$ can be defined inductively via
these relations.

\subsection{Definition of quantum Schubert polynomials}
\label{qsbdef}
We review the definition of quantum Schubert polynomials $\fS_w(x,q)$ 
given in \cite{fgp}.
Consider the same Dynkin diagram as in section \ref{usbdef}, 
with vertices labeled, $x_1,\ldots,x_n$.  Let $q_i$ be the (degree $2$) path
covering $x_i$ and $x_{i+1}$.  Then define $e^q_k(l)$ to be the sums of
monomials in $x$ and $q$ covering exactly $k$ of the $x_1,\ldots,x_l$.
For $w\in S_n$, define $\fS_w(x,q)$ to be the polynomial resulting from 
the substitutions $e^q_k(l)$ for $e_k(l)$ into the decomposition
of the ordinary Schubert polynomial $\fS_w(x)$ given by (\ref{spsum}).

It is easy to see from these descriptions that
$e^q_k(l)$ and $\fS_w(x,q)$ are the
 polynomials resulting from the substitution $g_i[1]=q_i$ and
$g_i[j]=0$ for $j\geq 2$ into $E_k(l)(g)$ and $\fS_w(g)$, respectively,
and that $\fS_w(x,0)=\fS_w(x)$,
 so that the quantum Schubert polynomials specialize further
to the classical Schubert polynomials.

\subsection{A degeneracy locus formula} \label{deg}
The universal Schubert polynomials are the solution to the following
degeneracy locus problem.
\begin{thm}
\label{degen}
Let $X$ be a Cohen-Macaulay scheme.  Consider maps of vector bundles 
$$V_1\rarr V_2 \rarr\cdots\rarr V_n = E_n \rarr E_{n-1}\rarr\cdots\rarr E_1$$
on $X$ where each $V_i$ and $E_i$ is of rank $i$, and
the $V_i$ are trivial vector bundles.  For each $w\in S_{n+1}$,
there is a degeneracy locus
$$\Omega_w = \{x\in X: \mbox{rank}_x(V_p\rarr E_q)\leq r_w(q,p) \mbox{ for }
1\leq p,q\leq n\}$$
where $r_w(q,p)$ is the number of $i\leq q$ such that $w(i)\leq p$.
Let  $[\Omega_w]$ be the cohomology class associated to $\Omega_w$ with
the scheme structure given locally by vanishing of determinants.
Assume that $\Omega_w$ is of the expected codimension $l(w)$.
Then setting $c_k(l)=c_k(E_l)$ the $k$th chern class of $E_l$, 
we have $[\Omega_w]=\fS_w(c)$.
\end{thm}

This is a consequence of the general result, Proposition 3.1 in  \cite{fu},
because all Chern
classes of $V_i$ vanish for $1\leq i\leq n$.

We have the following bundle maps on $\Pone\times\R$:
\begin{equation} \label{dualuniv}
V_1\rarr V_2 \rarr\cdots\rarr V_n = V_n \rarr 
A_{n-1}^*\rarr\cdots\rarr A_1^*
\end{equation}
given by the dual sequence to (\ref{univ}), where we have chosen a
fixed flag $V_\bullet$, and have denoted by $V_i$ the
corresponding rank $i$ vector bundle.  

Let $t\in\Pone$.  We can apply Theorem \ref{degen} to the bundle maps on 
$\R\eqq t\times\R$:
$$V_1\rarr V_2 \rarr\cdots\rarr V_n = V\otimes\oh_{t\times\R} \rarr 
(A_{n-1})_t^*\rarr\cdots\rarr (A_1)_t^*$$
given by restricting all bundles to $t\times\R$.

For the purpose of notation, set $C_k(l)$ to be the $k$th Chern class of the 
bundle $(A_l)_t^*$ on $\R$:
\begin{equation} \label{Ckl}
C_k(l)=c_k((A_l)_t^*).
\end{equation}
Define $Q_i[j]$ to be the classes in $H^*(\R)$
defined inductively
in terms of the $C_k(l)$ by the relations:
\begin{equation}\label{Qrecursion}
C_k(l) = C_k(l-1)+\sum_{j=0}^{k-1} C_{k-j-1}(l-j-1)Q_{l-j}[j]+
Q_{l-k}[k].
\end{equation}
These relations between the classes $C_*(*)$ and $Q_*[*]$ 
are specializations of the recursive relations described for the variables
$c_*(*)$ and $g_*[*]$ in (\ref{recursion2}).

By the corollary to Theorem \ref{gpschubert}, we have 

\begin{prop} 
\label{univQ}  Fix any point $t\in\Pone$. For $1\leq j\leq n-1$, let
 $(A_j)_t$ be the restriction of 
the tautological subbundles on $\Pone\times\R$ to $t\times\R$,
with $Q_i[j]$ classes in $H^*(\R)$ 
defined in terms of the Chern classes of the bundles $(A_j^*)_t$
as in (\ref{Qrecursion}).
Then the class $\mu_w(\bd)=[\overline{\Omega}_w(t)]$ in $H^*(\R)$ is given by
the formula $$\mu_w=\mu_w(\bd)=\fS_w(c_k((A_l^*)_t))=\fS_w(Q).$$
\end{prop}

\subsection{A geometric interpretation of the variables $g_*(*)$}
In this section, we abuse notation by letting $g_i[j]$ denote
the specialization of the variable $g_i[j]$ to certainly cohomology
classes on a scheme $X$.  We consider the setting of flags of vector bundles
on $X$:
$$E_n \rarr E_{n-1}\rarr\cdots\rarr E_1,$$
and show that, in some sense, the $g_*(*)$
measure how far the maps of vector bundles $E_l\rarr E_k$
are from being surjective. We apply Proposition \ref{gij} of this section 
to the proof of Proposition \ref{j2} in section \ref{proof}.

We define some polynomials which we use in this section.
For $a<b$, define $E_i(a,b)(g)$ to be the sum of all monomials in paths
$g_*(*)$ covering exactly $i$ of the vertices $x_a,\ldots, x_b$, no vertex
more than once.  When the variables (or classes) $g$ are understood, we write
$E_i(a,b)$.  With this
notation, $E_i(j)=E_i(1,j)$ and $E_i(a,b)=0$ for $i>b-a+1$.

\begin{prop} \label{gij}
Let
$$E_n \rarr E_{n-1}\rarr\cdots\rarr E_1$$ be a sequence of vector bundles 
on a scheme $X$ with $\mbox{rank }E_i =i$.  Let $c_k(l)=c_k(E_l)$ and
let $g_i[j]$ be defined in terms of the $c_*(*)$ by the relations
(\ref{recursion2}).
If $E_l\rarr E_k$ is a surjection of vector bundles for some $l\geq k$, then
 $g_i[j]=0$ for all $i,j$ satisfying $i<k+1\leq i+j\leq l$.
\end{prop}
\bpf  We proceed by induction on $l$.  The base case $l=k$ holds trivially
We assume that the result holds for all $l'<l$. 
If $E_l\rarr E_k$ is surjective, then so is
$E_{l'}\rarr E_k$ for $k\leq l'<l$.  In particular, 
$\ker(E_{l'}\rarr E_k)$ is a vector bundle of rank $l'-k$.
Furthermore, by the induction hypothesis,
we know that $g_i[j]=0$ for $i<k+1\leq i+j <l$, so it suffices to prove 
the result for $i+j=l$.

We  prove and use the claim:
\begin{lem} \label{Ekl} Assume that the proposition holds for $l'<l$.  Then
$c_i(\ker(E_l\rarr E_k)) = E_i(k+1,l')(g)$ for $l'\leq l$,
where the polynomials $E_i(j)(g)$ are defined by (\ref{recursion}).
\ref{usbdef}.
\end{lem}
\bpf
Denote by $E(k+1,l')(g)$ the polynomial $\sum_{i=0}^{l'-k}E_i(k+1,l')(g)$
The $i$th degree component of $c(E_{l'})$ differs from the $i$th degree of
 the product $c(E_k)E(k+1,l')$ by exactly all degree $i$ monomials
in $g$ containing a path $g_i[j]$ with $i\leq k$
and $k+1\leq i+j\leq l'\leq l$, i.e. paths joining $x_i$ for $i\leq k$
to $x_{i+j}$ for $k+1\leq i+j$.  By assumption, these are all zero except
possibly when $i+j=l'=l$.  
But $g_i[j]$ is of degree $j+1=l-i+1\geq l-k+1$ since
$i\leq k$, so such $g_i[j]$ do not occur in $E(k+1,l)$.  Therefore
we have the equality $c(E_k)E(k+1,l')=c(E_{l'})$ for $l'\leq l$ as needed.
\epf

With this result, we can complete the proof of the proposition.  We need
to prove that $g_i[j]=0$ for $i\leq k, i+j=l$.  We use induction on $j$
for $j\geq l-k-1$.  Assume that the result holds for $j'<j$.
By the lemma, we know that $c(E_k)E(k+1,l)=c(E_l)$.
For any $j$ the difference between the degree $j+1$ parts is the sum
of all monomials of degree $j$ containing a variable $g_i[j']$ satisfying
$i<k+1\leq i+j'\leq l$.  The variable $g_i[j']$ may only occur
 when $j'\leq j$  since $g_i[j']$ has degree
$j'+1$. By the first induction hypothesis, $g_i[j']=0$ except when
 $i+j'=l$, and by the second inductive hypothesis, $g_{l-j'}[j']=0$
when $j'<j$.  Therefore, the only term remaining is $g_{l-j}[j]$, but
on the other hand we began with an equality, so the difference
$g_{l-j}[j]$ is forced to be zero as well.
This concludes the proof.
\epf

\begin{cor}  With the same notation and hypotheses of Proposition
\ref{gij},
$$c(\ker(E_l\rarr E_k)) = \sum_{i=0}^{l-k}E_i(k+1,l)(g).$$
\end{cor}

\section{The results}
\label{results}

In this section, we state and assume the Main Proposition, which is proved in
section \ref{mainprop}. We  use the degeneracy
locus formula applied to bundles on the hyperquot scheme as in 
Proposition \ref{univQ}. 

For any $w\in S_n$, let $\sigma_w$ be the corresponding Schubert class in
the small quantum cohomology ring $QH^*(\F(n))$.  Throughout this section,
consider the universal sequence of sheaves (\ref{dualuniv}).  Let
the $Q_i[j]$ be classes in $H^*(\R)$, defined in terms of the Chern classes
of $(A_i^*)_t$, for any fixed $t\in\Pone$ as in the statement of 
Proposition \ref{univQ}.
Let $\mu_w$ be the classes in the hyperquot schemes as described in 
section \ref{schubert}.

We have the following formulation of the definition of the quantum product,
given by the definition in (\ref{multmap}) and Corollary \ref{Omegaprod}
to Theorem \ref{gpschubert}.

\begin{equation}
\label{qprod}
\sigma_{w_1}*\ldots * \sigma_{w_N} =
\qsum (\mu_{w_1} \cdots \mu_{w_N}\cdot\mu_w)_\bd \sw.
\end{equation}
where $(\mu_{w_1} \cdots \mu_{w_N}\cdot\mu_w)_\bd$ is the intersection
product of the classes $\mu_*$ in the hyperquot scheme $\R$.

The proofs of the results rely on two simple lemmas.

\begin{lem}
\label{polynomial}
Let $f$ be a polynomial with integer coefficients in variables $\nu_w$
for $w\in S_n$. Then
$$f(\sigma_*) = \qsum (f(\mu_*)\cdot\mu_{w})_\bd \sw$$
in the quantum cohomology ring of $\flag$.
\end{lem}
\bpf
Writing $f$ as a sum of monomials, the result follows 
immediately from (\ref{qprod}) and the additivity of both the quantum cohomology ring and the 
cohomology of the hyperquot schemes.
\epf

The following lemma is proved by an argument of Bertram, cf. \cite{b}.
\begin{lem}
\label{id}
For $w\in S_n$, $\sigma_{w}=\qqsum (\mu_{w}\cdot\mu_{w'})_\bd \sww$
in $QH^*(\F(n))$.
\end{lem}
\bpf
When $\bd=0$, i.e. $d_i=0$ for all $i$, $\hquot=\F(n)$, and
the $\mu_w$ are the ordinary Schubert cycles on the flag manifold.
In this case, $\mu_{w}\cdot\mu_{w'}=1$ when $\dualww=w$ and zero otherwise.

Therefore, it is enough to show that $(\mu_{w}\cdot\mu_{w'})_\bd=0$
whenever $\bd\neq 0$.  Recall that $(\mu_{w}\cdot\mu_{w'})_\bd
=\langle\Omega_{w},\Omega_{w'}\rangle_\bd$.  This Gromov-Witten number is zero
since only two points in $\Pone$ have been fixedl
\epf

For the remainder of the paper, we fix some notation.
For any multidegree $\bd$, let the classes $Q_i[j]\in H^*(\R)$ be 
defined in terms of the Chern classes of the tautological bundles by 
(\ref{Qrecursion}), where the
multidegree $\bd$ of $\R$ is understood by the appearance of
the $Q$'s in intersection products $(\cdots)_\bd$ on the
corresponding hyperquot scheme $\R$.  
For any polynomial $P(g)$ in the variables $g_i[j]$ for
$j\geq 0, i\geq 1$ and $i+j\leq n$, denote by $P(Q)$ the corresponding
polynomial in the classes $Q_i[j]$.  Define $P(x,q)$ to be the polynomial 
resulting from
setting $g_i[0]=x_i, g_i[1]=q_i$, and $g_i[j]=0$ for $j\geq 2$ in $P(g)$.

With this notation, we have the following claim.

\vspace{+10pt}
\noindent {\bf Main Proposition.}
For any multidegree $\bd$, let the classes $Q_i[j]\in H^*(\R)$ be 
defined in terms of the Chern classes of the tautological bundles by 
(\ref{Qrecursion}).
Let $P(g)$ be any polynomial in $g_i[j]$.  Then for $w\in S_n$,
$$P(x,q) = \qsum (P(Q)\cdot\mu_w)_\bd \sw \mbox{ in }QH^*(\flag).$$

The Main Proposition is proved in the next section, with the exception of 
a key lemma, whose proof is found in section \ref{proof}.
We conclude this section by proving the presentation of 
$QH^*(\flag)$ and the quantum Giambelli formula
as two easy corollaries of the Main Proposition.

For any $1\leq i\leq j\leq n$ with $i+j\leq n$, let 
$E_k(l)(g)$ be the polynomial in the variables $g_i[j]$
 as in section \ref{usbdef}.
Recall that from our observations in section \ref{qsbdef},
 $e_k^q(l)=E_k(l)(x,q)$, and that $e_k(l)$ is the $k$th 
elementary symmetric polynomial in the variables $x_1,\ldots,x_l$.
\begin{thm}
The small quantum cohomology ring of the complete flag manifold
has a presentation:
$$QH^*(\flag)=\Z[x_1,\ldots,x_n,q_1,\ldots,q_{n-1}]/I_q$$
where $I_q=(e_1^q(n),\ldots,e_n^q(n))$.
\end{thm}
\bpf By \cite{st}, it suffices to produce $n$ relations, 
 which specialize to the $n$ relations in $H^*(\flag)$ upon
setting $q_1=\cdots = q_{n-1}=0$.  
Recall that the relations of $H^*(\flag)$
are given by $e_1(n),\ldots,e_n(n)$.  
Since $A_n^*=V_n^*$ is a trivial vector bundle, we know that
$c_i(A_n^*)=E_i(n)(Q)=0$ for every $i$.  For 
$1\leq i\leq n$, we apply the Main Proposition to the polynomial
$E_i(n)(g)$ to get
$$0 = \qsum (E_i(n)(Q)\cdot\mu_w)_\bd \sw = E_i(n)(x,q)=e_i^q(n).$$
It is clear that $E_i(n)(x,0)=e_i(n)$, so that $e_i^q(n)$ specializes
to $e_i(n)$. This concludes the proof of the theorem.
\epf

\begin{thm} In the small quantum cohomology ring of $\F(n)$, there
is the  quantum Giambelli formula,
$$\sigma_w = \mathfrak{S}_w(x,q).$$
\end{thm}
\bpf
By the degeneracy locus formula, $\mu_w= \mathfrak{S}_{w}(Q)$.  Putting this
into  Lemma \ref{id} and applying the Main Proposition
to the polynomial $\fS_w(g)$ gives:
\begin{eqnarray*}
\sigma_{w}
&=&\qqsum (\mathfrak{S}_{w}(Q)\cdot\mu_{w'})_\bd \sww \\
&=& \mathfrak{S}_w(x,q).
\end{eqnarray*}
\epf

\section{Main Proposition}
\label{mainprop}
The Main Proposition asserts that there is some sort of correspondence between
 classes $Q_i[j]$ in $H^*(\R)$ over various $\bd$ with elements
in $QH^*(\flag)$.
In particular, it appears that there should be a relationship between
$Q_i[0]$ and $Q_i[1]$ with the $x_i$ and the deformation variables $q_i$ 
in  $QH^*(\flag)$, respectively.
Moreover, the classes $Q_i[j]$ with $j\geq 2$ should give,
in some sense, zero contribution to quantum cohomology.  

In this
section, we obtain an understanding of some classes $Q_i[j]$ for 
a handful of multidegrees $\bd$, 
and state Proposition \ref{q1mult}, which computes a particular
type of intersection on the hyperquot scheme.  Using
the structure of the
quantum cohomology ring, the cohomology rings of hyperquot schemes, and
the correspondence between $\Z[c]$ and $\Z[g]$ in (\ref{recursion2}), we show
that this is enough to prove the Main Proposition.  

The general argument for the analog of Proposition \ref{q1mult} for
partial flag manifolds is postponed to section \ref{qcpf}.  
For complete flag manifolds, we have additional structure given by describing
the classes  $Q_i[1]\in H^*(\R)$ for some $\bd$ 
as classes of particular degeneracy loci, and
study them geometrically in sections \ref{other}
and \ref{proof}.  Proposition \ref{q1mult} follows as a 
special case of Proposition \ref{Qeqn} in section \ref{proof}

For any polynomial $P(g)$ in the variables $g_i[j]$ for
$j\geq 0, i\geq 1$ and $i+j\leq n$, write $P^k(g)$ 
for the resulting polynomial  after setting $g_i[j]=0$ for $j\geq k+1$ in
$P(g)$. Note that $P^k(g)$ is a polynomial in variables
$g_i[j]$ with $j\leq k$.  With this notation,  
$P^0(x,q)=P(x)$ and $P^k(x,q)=P(x,q)$ for $k\geq 1$.

In order to prove the Main Proposition, we prove an auxiliary result:
\begin{prop}
\label{Pk}  Let $P(g)$ be any polynomial in $g_i[j]$. For $k\geq 0$,
$$P^k(x,q)=\qsum (P^k(Q)\cdot\mu_w)_\bd \sw.$$
\end{prop}
The Main Proposition follows as the special case $k=n$ of Proposition \ref{Pk}.
We use induction on $k$ to prove Proposition \ref{Pk},
 using Lemma \ref{polynomial} in almost every step.
We prove the result for $k=0$,  and then use induction. 

For the remainder of the paper, we set $y_i=Q_i[0]$, a class in $\R$, and
define $|\bd| = \sum d_i$.

\begin{prop} 
\label{x}
If $P(g)$ is a polynomial in the variables $g_i[0]$, then 
$$P(x)=\qsum (P(y)\cdot\mu_w)_\bd \sw \mbox{ in }QH^*(\flag).$$
\end{prop}
\bpf
In each hyperquot scheme $\R$, 
we can write $y_i = Q_i[0]= \mu_{s_i}-\mu_{s_{i-1}}$.
There is no quantum deformation of divisor classes so that the
quantum Giambelli formula for $w=s_i$ holds, $\sigma_{s_i} = x_1 + \cdots + x_i$,
so that $x_i = \sigma_{s_i}-\sigma_{s_{i-1}}$.  Therefore,
it suffices to prove the statement for a polynomial in the $\sigma_{s_i}$'s.
The statement follows as a special case of 
Lemma \ref{polynomial} where each $w$ is some $s_i$.
\epf

\begin{prop}
\label{q1mult}  For $1\leq i\leq n-1$, $(Q_i[1]\cdot \mu_{w_0})_\bei =1$,
where $w_0$ is the permutation of longest length in $S_n$.
\end{prop}

It turns out that these $n-1$ intersections are the only ones that 
must be computed  in order to see the effect of the classes $Q_i[1]$ 
on quantum cohomology.  There is an analagous result for partial flag
manifolds in section \ref{qcpf}.

We prove some facts  regarding the classes $Q_i[j]$ for $j\geq 2$
in $H^*(\R)$.  
\begin{lem} \label{j2lem}
If $d_k=0$ for some $k$ then $Q_i[j]=0$ when $i$ and $j$
satisfy:
\begin{enumerate}
\item $i\leq k$ and
\item $i+j\geq k+1$.
\end{enumerate}
\end{lem}
\bpf
 If $d_k=0$, then $A_k\rarr V_n=A_{n}$ is a vector bundle
inclusion, or equivalently,  $A_n^*\rarr A_k^*$ is a surjection.
Then a direct application of Proposition \ref{gij} gives the result.
\epf

\begin{prop}
\label{j2}
\begin{enumerate}
\item
If $|\bd| \leq j-1$, then $Q_i[j]=0$ in $H^*(\R)$ for 
every $1\leq i\leq n-j$. 
\item
If $j\geq 2$, then in $QH^*(\flag)$
$$ \qsum (Q_i[j]\cdot \mu_w)_\bd \sw    =0.$$
\end{enumerate}
\end{prop}

\bpf
Part 1 of the proposition follows from the fact 
that if $Q_i[j]\neq 0$ for some $i, j$, then by Lemma \ref{j2lem}, $d_k\neq0$
for $k=i,\ldots,i+j-1$.  In particular, $d_k\geq 1$ so that
$|\bd|\geq \sum_{k=i}^{i+j-1} d_k \geq \sum_{k=i}^{i+j-1} 1\geq j$.

To prove Part 2, we consider two cases:
\begin{enumerate}
\item $|\bd| \leq j-1$.  In this case, by part 1, $Q_i[j]=0$ so that
$(Q_i[j]\cdot \mu_w)_\bd=0$ for all $w$ and $\bd$.
\item $|\bd| \geq j$.  
The product $(Q_i[j]\cdot \mu_w)_\bd$ can only
be nonzero in the correct dimension, when $(j+1)+l(w)=\dim(\R)
=\binom{n}{2}+2|\bd|$. Since $l(w)\leq \binom{n}{2}$ for all $w$,
and $j\geq 2$ implies that $2j>j+1$, this
situation can never occur, so that the intersection is zero by
dimension considerations.
\end{enumerate}
This concludes the proof of the proposition.
\epf

Part 1 of Proposition \ref{j2} is the only direct knowledge of the geometry of
the classes $Q_i[j]$ for $j\geq 2$ needed for the results.  Besides the proof
of Proposition \ref{q1mult}, which requires a detailed geometric description of
the classes $Q_i[1]$ and $\mu_w$, the remainder
of the proof of the Main Proposition is given by algebraic manipulations on the
classes $\mu_w=\fS_w(Q)$ and $Q_i[j]$ in $H^*(\R)$ for all $\bd$, 
and properties of the quantum cohomology ring.  

\begin{lem}
\label{muk+1}
Let $k$ be an integer such that Proposition \ref{Pk} holds for $k$.   Then
$\sigma_w = \fS_w(x,q)$ in $QH^*(\flag)$ for $l(w)\leq k+2$.
\end{lem}

\begin{lem}
\label{Qisigma}   Let $k$ be an integer such that  Proposition \ref{Pk} 
holds for $k$. Then the class $Q_i[k+1]$ can be written
as a polynomial in those classes $\mu_w$ with $l(w)\leq k+2$.
\end{lem}

We first prove Lemma \ref{muk+1}.  For
$l(w)\leq k+2$, the polynomial $\fS_w(Q)$ involves classes
$Q_i[j]$ where $j\leq k+1$.  Therefore, 
we can write $\fS_w(Q)=\fS_w^k(Q)+\sum_{i\in I}
Q_{i}[k+1]$ for some sequence $\{i\in I\}$, possibly empty, with the 
$i_j\in\{1,\ldots,n-k-1\}$ not necessarily
distinct.  Then the additivity of the cohomology of the hyperquot scheme,
and the fact that $\mu_w=\fS_w(Q)$ gives:
\begin{eqnarray*}
\sigma_w&=& \qqsum (\fS_w(Q)\cdot\mu_{w'})_\bd \sww\\
&=& \qqsum (\fS^k_w(Q)\cdot\mu_{w'})_\bd \sww
+ \sum_{i_j} \qqsum (Q_{i}[k+1]\cdot\mu_{w'})_\bd \sww\\
&=& \fS_w^k(x,q)+\sum_{(i,k+1)=(i,1)} q_{i} \\
&=& \fS_w(x,q).
\end{eqnarray*}
The third equality is an immediate application of the property that $k$ 
satisfies the statement of Proposition \ref{Pk}, Proposition \ref{Qeqn}, and
Part 2 of Proposition \ref{j2}.  The final equality follows from the definition
of the polynomials $\fS^k_w(x,q)$ and $\fS_w(x,q)$.
\epf

We now prove Lemma \ref{Qisigma}.  Let $\beta_{i,k}\in S_n$
be the permutation with $\mu_{\beta_{i,k}}=C_{k+2}(i+k+1)$, where the
classes $C_*(*)$ are defined by (\ref{Ckl}). 
In cycle notation $\beta_{i,k}=(i \, i+1 \cdots k+i+2)$, and 
$l(\beta_{i,k})=k+2$. In the notation of \cite{cf1}, 
$\beta_{i,k}=\alpha_{k+3,k+i+2}$.

By the recursion in (\ref{Qrecursion}) and the fact that $\mu_{\beta_{i,k}}=
\fS_{\beta_{i,k}}(Q)$, we
can write $$Q_1[k+1]+\cdots +Q_i[k+1]=\mu_{\beta_{i,k}}-
\fS_{\beta_{i,k}}^k(Q).$$ 
By the induction hypothesis,
every  $Q_i[j]$ with $j\leq k$ can be written
as a polynomial in $\mu_w$ with $l(w)\leq k+1$.  Therefore, since
$\fS_{\beta_{i,k}}^k(Q)$ is by definition a polynomial in
$Q_i[j]$ with $j\leq k$, it  can be as well.
Therefore, by induction on $i$, we see that each $Q_i[k+1]$ can be written
as a polynomial in $\mu_w$ with $l(w)\leq k+2$.  
\epf

We are ready to complete the inductive
proof of  Proposition \ref{Pk}.  The base case
$k=0$ was proved in Proposition \ref{x}.  Let $k$ be
such that Proposition \ref{Pk} holds.  By Lemma \ref{Qisigma}, every $Q_i[j]$
with $j\leq k+1$ can be written as a polynomial in $\mu_w$ with
$l(w)\leq k+2$.   By definition,
 $P^{k+1}(Q)$ is a polynomial in these  $Q_i[j]$,
so that it can be written as $P^{k+1}(Q)=\hat{p}(\mu_*)$ with $l(w)\leq
k+2$.

We apply Lemma \ref{polynomial} to $P^{k+1}(g)$ to get
$$\qsum (P^{k+1}(Q)\cdot\mu_w)_\bd \sw
= \qsum (\hat{p}(\mu_*)\cdot\mu_w)_\bd \sw
= \hat{p}(\sigma_*)=P(x,q)$$
where the final equality follows from Lemma \ref{muk+1} and
the fact that $\hat{p}(\sigma_*)$ is a polynomial in $\sigma_w$ with
$l(w)\leq k+2$.
\epf

Except for the proof  of Proposition \ref{q1mult}, 
this concludes the proof of Proposition \ref{Pk} and the
Main Proposition, and hence of the presentation
of $QH^*(\flag)$ and of the quantum Giambelli formula.

\section{More degeneracy loci on the hyperquot scheme}
\label{other}

In this section, we define and study
certain degeneracy loci on $\R$ which allow us
to geometrically and explicitly understand the classes $Q_i[1]$ in $H^*(\R)$
for the various hyperquot schemes.  These loci are crucial to the proof
of Proposition \ref{q1mult}, which is what remains to be proved.

Over $\Pone \times \hquot$, there is the universal sequence of sheaves
(\ref{univ}).
For $1\leq i\leq n-1$, define $W_i \subset \Pone \times \R$
to be locus 
over which $\mbox{rank} (A_i \rarr A_{i+1}) \leq {i-1}$, with the
natural scheme structure is given by the vanishing of determinants.
Then for $p\in\Pone$, define $W_i(p)=W_i \cap (p\times \R)$.
Identifying $p\times \R$ with $\R$, we view $W_i(p)$ as a subscheme of $\R$.
We use points $p$ for the loci $W_i(p)$, while we use points $t$ for the loci
$\overline{\Omega}_{w}(t)$ introduced in section \ref{schubert}.

We prove a general position
result in the spirit of the argument of section \ref{schubert},
which shows that $W_i(p)$ is of expected complex codimension two, thus
giving a class in the cohomology of the hyperquot scheme.  Furthermore,
we show that this class is independent of the choice
of the point $p$, and that this class is given by $Q_i[1]$.

The following lemma is the analog of Lemma \ref{hO} to these degeneracy
loci.
\begin{lem}
\label{hq}
Let $\be$ be a multiindex, with $h_\be$ and $\pi$ as
defined in section \ref{boundary}.  Then
$$h_\be^{-1}(W_i(p))= \pi^{-1}(\Pone\times W_i(p)) 
\cup \widetilde{W}_i(p)$$
with $\widetilde{W}_i(p)$ the degeneracy locus inside 
$\pi^{-1}(p\times \Rde)$
given by $$\mbox{rank }(K_i\rarr K_{i+1})\leq i-1.$$
This equality is 
scheme-theoretic away from the intersection.
\end{lem}
\bpf
We recall the construction of the map $h_\be$
as stated in section \ref{hyperquot}.
For $1\leq i\leq n-1$, we have
$(1\times h_\be)^*A_i^\bd = \widetilde{A}_i^\be$,
so that $$h_\be^{-1}(W_i(p))=\{ y\in\U_\be | \mbox{rank}_{(p,y)}
\widetilde{A}_i^\be \rarr\widetilde{A}_{i+1}^\be \leq i-1 \}.$$
Outside $\widetilde{\Delta}$, for $j=i,i+1$,
$\psi^*A_j^\bde$ is isomorphic to 
$\widetilde{A}_{j}^\be$, while over the locus $\widetilde{\Delta}$, we have
 $\widetilde{A}_{j}^\be=K_j$.  These two observations give the lemma.
\epf

In fact we can say more:
\begin{cor}
\label{Wei} For a multiindex $\be=(e_1,\ldots,e_{n-1})$,
\begin{enumerate}
\item
If $e_i=0$, then $h_\be^{-1}(W_i(p))= \pi^{-1}(\Pone\times W_i(p))$
as schemes. 
\item
If $e_i>0$, then $h_\be^{-1}(W_i(p))= \pi^{-1}(\Pone\times W_i(p)) 
\cup \pi^{-1}(p\times \hquote)$ as schemes away from the intersection.
\end{enumerate}
\end{cor}
\bpf
Recall that $\widetilde{W}_i(p)$ is defined as the locus
in $\pi^{-1}(p\times \hquote)$
given by the condition
$\mbox{rank }(K_i\rarr K_{i+1})\leq i-1.$  In particular,
we see that elements in $\widetilde{W}_i(p)$ satisfy the condition
$\mbox{rank }(K_i\rarr V^*)\leq i-1.$

Further recall that $\U_\be(p)= \pi^{-1}(p \times \hquote)$
is defined as the locus where $\mbox{rank }(K_i\rarr V^*) =  i-e_i$
for $1\leq i\leq n-1$.

By these definitions, we see that when $e_i=0$,
$\widetilde{W}_{i}(p)$ is empty, giving the first part of the claim.

The locus where  $\mbox{rank }(K_i\rarr V^*) =  i-e_i$ is clearly
contained in the locus where $\mbox{rank }(K_i\rarr K_{i+1})\leq i-1$.
Therefore, $\U_\be(p))= \pi^{-1}(p \times \hquote)$ 
is contained in $\widetilde{W}_i(p)$.
This gives the second part of the claim.
\epf

In order to intersect these loci, we  need the following theorem 
about general position, extending the general position results
of \cite{cf1} as stated in Theorem \ref{gpschubert}.  

\begin{thm} (Moving lemma)
\label{moving}
For $i_1,\ldots,i_M$ in $\{1,\ldots,n\}$,
permutations $w_1,\ldots,w_N\in S_n$, general translates of 
$\Omega_{w_k}\subset \F(n)$, and distinct points 
$t_1,\ldots,t_N,p_1,\ldots,p_M\in \Pone$, the
intersection 
\begin{equation}
\label{intersection}
\bigcap_{j=1}^M W_{i_j}(p_j) \cap 
\bigcap_{k=1}^N \overline{\Omega}_{w_k}(t_k)
\end{equation} is either empty or
has pure codimension $2M+\sum_{k=1}^N l(w_k)$.  Here the
$i_j$ and the $w_k$ are not necessarily distinct.

\end{thm}
\bpf
For $M=0$, this is Theorem \ref{gpschubert}.  We only need
to consider the case where $M\geq 1$.
The proof is by induction on $\bd=(\dd)$.  

The base case is when $d_1=\cdots=d_{n-1}=0$, so that
 $\hquot = \F(n)$.  In particular,
$A_i\rarr A_{i+1}$ is a vector bundle inclusion for each $i$.
Thus $W_i=\emptyset$ for all $i$, and hence $W_{i_1}(p_1)$ (which appears
in the intersection if $M\geq 1$) is empty so
that the entire intersection is empty.  

For two multiindices,  write  $\mathbf{f}<\bd$ when 
$f_j\leq d_j$ for every $1\leq j\leq n-1$
and $f_k<d_k$ for some $1\leq k\leq n-1$.  Assume that the result
holds for all such $\mathbf{f}$.

Since any $W_i$ is the locus where $A_i\rarr A_{i+1}$ drops rank,
we must have $\mbox{rank }(A_i\rarr V)\leq i-1$ over $W_i$ as well.
In particular, it follows from the second part of Theorem \ref{Ue} 
that for any $1\leq i\leq n-1$ and $p\in\Pone$, $W_i(p)$ is contained
in the boundary of the hyperquot scheme.

Let $L=\sum_{k=1}^N l(w_k)$. 
By assumption, $M\geq 1$  so
that the intersection (\ref{intersection}) is also contained in the boundary.
Since $\bigcup_\be h_\be(\U_\be)$ covers the boundary,
it suffices to show that for every $0<\be<\bd$,
$$\mbox{codim}_{\R}
\left(h_\be(\mathcal{U}_\be) \cap \bigcap_{j=1}^M W_{i_j}(p_j) \cap 
\bigcap_{k=1}^N \overline{\Omega}_{w_k}(t_k)\right) \geq 2M+L.$$
The map $h_\be$ is birational onto its image, so we only need to show,
for every $\be$, that
\begin{eqnarray*}
\lefteqn{\mbox{codim}_{\mathcal{U}_\be}
\left(\bigcap_{j=1}^M h_\be^{-1}(W_{i_j}(p_j)) \cap 
\bigcap_{k=1}^N h_\be^{-1}(\overline{\Omega}_{w_k}(t_k))\right) }\\
&&\geq 2M+L-(\dim\R -\dim\U_\be)\\
&&=2M+L+1-\sum e_i(1+e_i-e_{i-1}).
\end{eqnarray*}

$\widetilde{W}_{i_j}(p_j)$ and $\widetilde{\Omega}_{w_k}(t_k)$ 
are supported on  $p_j\times\hquote$ and $t_k\times\hquote$, respectively,
By  Lemma \ref{hO}, Lemma \ref{hq}, and the fact that
 $p_1,\ldots,p_M,t_1,\ldots,t_N$ are all distinct points in $\Pone$,
we see that after a possible renumbering,
 the only possible nonempty intersections are of the type
\begin{equation}
\label{none}
\pi^{-1}\left(\Pone\times \bigcap_{j=1}^M W_{i_j}(p_j)\cap \bigcap_{k=1}^N
\overline{\Omega}_{w_k}(t_k)\right),
\end{equation}
of type
\begin{equation}
\label{O_N}
\pi^{-1}\left(t_N \times \bigcap_{j=1}^{M}W_{i_j}(p_j)) \cap 
\bigcap_{k=1}^{N-1} \overline{\Omega}_{w_k}(t_k)\right) 
\cap \widetilde{\Omega}_{w_N}(t_N),
\end{equation}
or of type
\begin{equation}
\label{Q_M}
\pi^{-1}\left(p_M \times \bigcap_{j=1}^{M-1}  W_{i_j}(p_j)) \cap
\bigcap_{k=1}^{N} \overline{\Omega}_{w_k}(t_k)\right) 
\cap \widetilde{W}_{i_M}(p_M).
\end{equation}

Since $\bd-\be<\bd$, by the induction hypothesis and the fact that 
$\pi$ is smooth,
intersections of type (\ref{none}) are of codimension $\geq 2M+L$ in
$\U_\be$. 

By Lemma \ref{hq} and Lemma \ref{codimO} and the induction hypothesis, 
intersections of type (\ref{O_N}) are of codimension $\geq 1+ 2M + L -l(w_N)
+l(\widetilde{w}_N^\be)$ in $\U_\be$.  But we know that
$l(w_N)-l(\widetilde{w}_N^\be)\leq \sum e_i$, so that such intersections 
are codimension
at least $2M+L-\sum e_i+1$ in $\U_\be$. 
 
By Lemma \ref{hq}, intersections of type (\ref{Q_M}) are
empty unless $e_{i_M}>0$, in which case $i=i_M$. 
From Corollary \ref{Wei}, we have 
$\widetilde{W}_{i_M}(p_M)=\pi^{-1}(p_M \times \hquote)$.
Therefore, (\ref{Q_M}) is either empty or
codimension $2(M-1)+L$ in $\U_\be(t_M)$, and hence codimension
$2(M-1)+L+1=2M+L-1$ in  $\U_\be$.

By part $1$ of Lemma \ref{es}, for any $\be$,
we have the inequalities 
$$2M+L \geq 2M+L-\sum e_i +1> 2M+L+1-\sum e_i(1+e_i-e_{i-1})$$
so that intersections of types (\ref{none}) and (\ref{O_N}) are empty
and  $2M+L-1 \geq 2M+L+1- \sum e_i(1+e_i-e_{i-1})$.  By part $2$ of
Lemma \ref{es},
this is an equality if and only if
 $|\be|=\sum e_i=1$, so that $\be=\bei$ for some $i$.  Therefore,
in the case $|\be|\geq 2$, intersections of type (\ref{Q_M}) are also empty,
and when $\be=\bei$,  we have $2M+L-1= 2M+L-(\dim\R -\dim \U_\bei)$,
which gives the needed codimension estimate.
\epf

In the course of the proof, we actually showed something stronger, that
the only nonempty intersection arise when $\be=\bei$, from type (\ref{Q_M}).
In particular, we have proven:

\begin{cor} \label{heimoving}
Consider the same hypotheses as in Theorem \ref{moving}, with
$\bei$ the multiindex with all zeros except a $1$ at the $i$th position.
Then
\begin{eqnarray*}
\lefteqn{h_\bei^{-1}\left(
\bigcap_{j=1}^M W_{i_j}(p_j) \cap 
\bigcap_{k=1}^N \overline{\Omega}_{w_k}(t_k)\right) } \\
& &= \bigcup_{i_j=i} \left(
\bigcap_{l \neq j}  \pi^{-1}(p_l\times W_{i_l}(p_l)) \cap
\bigcap_{k=1}^{N} 
 \pi^{-1}(\Pone\times \overline{\Omega}_{w_k}(t_k))\right) .
\end{eqnarray*}
\end{cor}

\begin{cor}
The class $[W_i(p)]\in H^4(\hquot)$ is independent of the choice of
the point $p\in\Pone$.
\end{cor}
\bpf  By definition, the  $W_i(p)$ are fibers of the morphism
 $W_i\subset \Pone\times\R\rarr \Pone$.   By Theorem \ref{moving}, 
the fibers are of complex codimension two. $W_i\rarr\Pone$ is in fact a fiber
bundle since the automorphism group of $\Pone$ is transitive.  
Therefore $[W_i(p)]$ is independent of the point $p$.
\epf

\begin{cor}
\label{Wirred}
For any $p\in\Pone$ and $1\leq i\leq n-1$, 
$W_i(p)$ is an irreducible scheme, with
$h_\be(\U_\bei(p))\subset W_i(p)$ an open subscheme.
\end{cor}
\bpf
We use the second part of Corollary \ref{Wei} applied to $\be=\bei$:
$$h_\be^{-1}W_i(p)=\pi^{-1}(\Pone\times W_i(p)) 
\cup \pi^{-1}(p\times \hquotei).$$

Since $\U_\bei(p) = \pi^{-1}(p \times \hquotei)$ is the preimage via the
smooth map $\pi$ of an irreducible scheme, it is also irreducible, as
its image $h_\bei(\U_\bei(p))$.
We observe that  $\codim_{\U_\bei} \U_\bei(p)=1$
and  $\codim_{\U_\bei} \pi^{-1}(\Pone\times W_i(p))  = 2$.
By Corollary \ref{Wei} $h_\bei^{-1}(W_i(p))$ is the union of these
two subschemes of $\U_\bei$, so that $W_i(p)$  is irreducible. 

Observe that
\begin{eqnarray*}
\dim W_i(p)= \dim\hquot -2 &=& \dim \Pone\times\hquotei -1 \\
&=&\dim\U_\bei -1 =\dim\U_\bei(p),
\end{eqnarray*} 

We know that $h_\bei $ maps $\U_\bei(p)$ isomorphically onto its image.
By the definition of $\U_\bei(p)$, it is a locally closed scheme, so that
its image $h_\bei(\U_\bei(p))$ is a union of locally closed schemes.
But $W_i(p)$ is irreducible, so has only one component, and therefore
$h_\bei(\U_\bei(p))$ is locally closed in $W_i(p)$.
Since $\dim\U_\bei(p)=\dim W_i(p)$, the image $h_\bei(\U_\bei(p))$ is
an open subscheme of $W_i(p)$.
\epf

\begin{prop}
\label{WQ}
Let the classes $Q_i[j]\in H^*(\R)$ be defined in terms of the Chern
classes of the tautological bundles $(A_l)_p$ on $\R$ by (\ref{Qrecursion}).
Then there is an equality $[W_i(p)]=Q_i[1]$.
\end{prop}
\bpf
By Theorem \ref{moving} and its corollaries, $[W_i(p)]$ is a class independent
of the point $p$ in $H^4(\R)$.  It is a degeneracy locus of expected
codimension two.  Since its expected dimension matches its actual dimension,
we can apply Porteous' formula so that $[W_i(p)]=c_2((A_{i+1})_p/(A_i)_p)$ where
$(A_i)_p$ is the bundle $A_i$ restricted to $p\times\R$.

By the definition of the classes $Q$ as in Proposition
\ref{univQ}, for any $1\leq l\leq n$, the Chern polynomial can be written
 $c(A_p^{l})=\sum_{k=0}^l (-1)^k E_k(l)(Q)$.
Let $b_1=c_1((A_{i+1})_p/(A_i)_p)$ 
and $b_2=c_2((A_{i+1})_p/(A_i)_p)$.  
Then $$b_1=c_1((A_{i+1})_p)-c_1((A_i)_p)=Q_{i+1}[0]$$
and the equation
$$c((A_{i+1})_p)=c((A_{i+1})_p/(A_i)_p)c((A_i)_p)$$
implies that  
\begin{eqnarray*}
E_2(i+1)(Q)=c_2((A_{i+1})_p)&=&c_2((A_i)_p)+(c_1((A_i)_p))b_1+b_2 \\
&=&E_2(i)+(-E_1(i))(-Q_{i+1}[0])+b_2. 
\end{eqnarray*}
We also have the recursion in equation (\ref{recursion})
$$E_2(i+1)(Q)=E_2(i)+E_1(i)Q_{i+1}[0]+Q_i[1]$$
so that $b_2=Q_i[1]$.
\epf

\section{Proof of Proposition \ref{q1mult}}
\label{proof}

In this section, we prove a more general statement than Proposition 
\ref{q1mult},
which gives a correspondence between the classes 
$Q_i[1]$ in the hyperquot schemes
$\R$ and the deformation variables $q_i$ in $QH^*(\flag)$.

For any multiindex $\bc=(c_1,\ldots,c_{n-1})$, write
$${\mathbf Q[1]}^\bc=Q_1[1]^{c_1}Q_2[1]^{c_2}\cdots Q_{n-1}[1]^{c_{n-1}}.$$
We prove the following proposition.  

\begin{prop}
\label{Qeqn} For any $\bc<\bd$,
$$ (\BQ^\bc \cdot \mu_{w_1}\cdots \mu_{w_N})_\bd = 
( \mu_{w_1}\cdots \mu_{w_N})_{\bd-\bc}.$$
\end{prop}

The proof is geometric in nature, involving degeneracy loci of types
$\overline{\Omega}_w(t)$ and $W_i(p)$ on 
$\hquot$ and an analysis of the boundary.  We use induction 
on $\bc$.
For $\bc=0$, this is Lemma \ref{x}, so we may assume that $c_m>0$ for some $m$.

Proposition \ref{WQ} states that $[W_i(p)]=Q_i[1]$ in $\R$ 
for any choice of 
point $p\in\Pone$.  By the moving lemma, we know that in order to intersect
these $W_i(p)$ it suffices to choose distinct points.
Choose $p_{i,j}$ and $t_k$ to be distinct points of $\Pone$ for 
$1\leq i\leq n-1, 1\leq j\leq c_i$, and $1\leq k\leq N$. 
Define $$Y:=\bigcap_{(i,j)\neq (m,1)} W_i(p_{i,j})\cap
\bigcap_{k=1}^N \overline{\Omega}_{w_k}(t_k)$$
where $m$ is such that $c_m>0$.

By the corollary to Theorem \ref{moving}, the left hand side of Proposition \ref{Qeqn} can be interpreted as the the intersection number
$[Y]\cdot [W_m(p_{m,1})]$ on $H^*(\R)$.

The idea of the proof
is to compute this intersection on a particular open subscheme of 
$W_m(p_{m,1})\subset\R$.  For the inductive step, 
we show that this intersection can 
be further computed on a smaller hyperquot scheme,  In order to do this, we
use the morphism  $h_\bem: \U_\bem \rarr \R$ and the open immersion
$\pi: \U_\bem\rarr \hquotem$.

By Lemma \ref{hq},
$h_{\bem}^{-1}(Y\cap W_m(p_{m,1}))$ 
is supported on $\pi^{-1}(p_{m,j}\times\hquotem)$ for
$1\leq j\leq c_i$.  Set-theoretically, we have 
$$\pi^{-1}(p_{m,j}\times W_m(p_{m,j'}))=
\pi^{-1}(p_{m,j'}\times W_m(p_{m,j}))$$ for any $j,j'$.  Hence,
$h_{\bem}^{-1}(Y\cap W_m(p_{m,1}))$ 
is supported on $\U_\bem(p_{m,1})=\pi^{-1}(p_{m,1}\times\hquotem)$ so that
the set-theoretic intersection 
$Y\cap W_m(p_{m,1})$ is contained in $h_\bem(\U_\bem(p_{m,1}))$.

By Corollary \ref{Wirred}, we know that  
$h_\bem(\U_\bem(p_{m,1}))$ is an  open subscheme of $W_m(p_{m,1})$.
Since (set-theoretically) $Y\cap W_m(p_{m,1})\subset 
h_\bem(\U_\bem(p_{m,1}))$, we have the equality of cycle 
intersections on $\R$
\begin{eqnarray*}
[Y]\cdot[W_m(p_{m,1})]&=&[Y]\cdot(h_\bem)_*[\U_\bem(p_{m,1})]\\
&=&  [(h_\bem)_*[(h_\bem)^{-1}(Y)\cap \U_\bem(p_{m,1})]
\end{eqnarray*}
where the second equality comes from the projection formula.

We know that $h_{\bem}$ maps $\U_\bem(p_{m,1})=\pi^{-1}(p_{m,1}
\times \hquotem)$ isomorphically onto its image, so that the right hand
quantity is equal to the intersection number 
$$ [(h_\bem)^{-1}(Y)\cap \U_\bem(p_{m,1})]\mbox{ on }
\U_\bem.$$

By Corollary \ref{heimoving} to
 Theorem \ref{moving}, we see that this last intersection is
\begin{eqnarray*}
\lefteqn{\bigcap_{(i,j)\neq(m,1)} \pi^{-1}(p_{m,1}\times W_i(p_{i,j}))
\cap \bigcap_{k=1}^N \pi^{-1}(p_{m,1}\times\overline{\Omega}_{w_k}(t_k)) }\\
& & =\pi^{-1}\left(p_{m,1}\times \bigcap_{(i,j)\neq(m,1)}W_i(p_{i,j})
\cap \bigcap_{k=1}^N\overline{\Omega}_{w_k}(t_k)\right)
\end{eqnarray*}
We claim that this is a scheme-theoretic equality.
By Lemma \ref{hO} and Lemma \ref{hq}, it suffices to show that this does not
intersect any of  $\pi^{-1}(t_k \times \overline{\Omega}_w(t_k))$ or
 $\pi^{-1}(p_{i,j}\times W_i(p_{i,j}))$.  Since all of
the points $t_k$ and $p_{i,j}$
are distinct, the only case to check is that of $p_{m,1}$.  But here, we
see that $\pi^{-1}(p_{m,1}\times (Y\cap W_m(p_{m,1})))$ 
is empty by the codimension results of
Theorem \ref{moving} applied to $Y\cap W_m(p_{m,1})$ on $\hquotem$.

All of the intersection points lie in the image of the
open immersion $\pi:\U_{\bem}\rarr \Pone\times\hquotem$.  Therefore,
via the  identification
$p_{m,1}\times\hquotem\eqq\hquotem$, 
$[Y]\cdot [W_m(p_{m,1})]$  is the length of the zero-dimensional
subscheme of $\hquotem$ given by
$$\bigcap_{(i,j)\neq(m,1)}W_i(p_{i,j})
\cap \bigcap_{k=1}^N\overline{\Omega}_{w_k}(t_k).$$
Alternatively, this is the intersection number 
$$\prod_{(i,j)\neq(m,1)}[W_i(p_{i,j})]
\cdot \prod_{k=1}^N[\overline{\Omega}_{w_k}(t_k)]
= \BQ^{\bc-\bem}\cdot
\mu_{w_1}\cdots\mu_{w_N} \mbox{ on } H^*(\hquotem).$$
Therefore, we have shown that
$$(\BQ^{\bc}\cdot\mu_{w_1}\cdots\mu_{w_N})_\bd = (\BQ^{\bc-\bem}\cdot
\mu_{w_1}\cdots\mu_{w_N})_{\bd-\bem}.$$  
Since $\bc-\bem < \bc$, the induction hypothesis on $\bc$ implies the result.
\epf

\section{Quantum cohomology of partial flag manifolds}
\label{qcpf}
In this section, we give the necessary ingredients to extend the arguments
to the quantum cohomology ring of partial flag manifolds.
 Let $N$ be the set $\{1\leq n_1<\ldots,n_m<n_{m+1}=n\}$.
Let $\pflag$ denote the partial flag variety corresponding
to flags of the form:
$$V_1 \subset V_2 \subset ... \subset V_m \subset V=\com^n$$
with dim $V_{m+1-i}=n-n_i$.   There is a universal sequence
of vector bundles on $\pflag$:
$$V_{\pflag}\rrarr Q_{m}\rrarr\cdots\rrarr Q_{1}$$
with $\mbox{rank }Q_j=n_j$, and each $Q_j\rarr Q_{j-1}$ a surjection, where
$Q_0:=0$ and $Q_{m+1}:=V_{\pflag}$.

\subsection{Classical cohomology of $\pflag$}
We review the ordinary cohomology of $\pflag$.  Let 
$S^{(N)}=\{w\in S_n:w(i)<w(i+1) \mbox{ for }i\not\in N\}$.
For $1\leq l\leq m+1$, let $x_{n_{j-1}+1},\ldots,x_{n_j}$ be the Chern roots
of the bundle $F_j:=\ker(Q_j\rarr Q_{j-1})$, and let 
$\sigma_i^j$ be the $i$th Chern
class of $F_j$ for $1\leq i\leq n_j-n_{j-1}$.    Then
 $e_k(n_l)$, the
$k$th elementary symmetric polynomial in $x_1,\ldots, x_{n_l}$, 
is symmetric in  $x_{n_{j-1}+1},\ldots,x_{n_j}$ for
every $1\leq j\leq l$, and can thus be written as a polynomial in the $\sigma_i^j$, which
 we denote by $\tilde{e}_k(l)(\sigma)$ or $\tilde{e}_k(l)$.

Let the polynomials $\fS_w(x)$ be as defined in section \ref{flag}.
Then it is a fact that for $w\in S^{(N)}$,
$\fS_w(x)$ can be written as a polynomial in $\sigma_i^j$, which
we write $\fS_w(\sigma)$.
For $w\in S^{(N)}$, consider the degeneracy locus
\begin{multline*}
\Omega_w(V_\bullet)=\{U_\bullet\in\pflag: \mbox{rank}_{U_\bullet}(V_p\otimes 
\oh_{\pflag}\rarr Q_q)\leq r_w(q,p) \\
\mbox{ for }1\leq p \leq n, q\in N\}.
\end{multline*}
This is a codimension $l(w)$ subvariety in $\pflag$ whose class is independent
of choice of flag, and is denoted by $[\Omega_w]$. Let $w^0\in S^{(N)}$
be the element of longest length, and let $\dualw=w^0w$.
We have the following classical results.
\begin{thm} The classes $\{[\Omega_w]\}_{w\in S^{(N)}}$ form an additive basis for
$H^*(\F)$.  Furthermore, for $w\in S^{(N)}$, the Schubert classes
$[\Omega_w]$ and $[\Omega_{\dualw}]$ are Poincar\'{e} dual.
\end{thm}

\begin{thm} $H^*(\pflag,\Z)\eqq \Z[\sigma_i^j]/(\e_1(m+1),\ldots,\e_n(m+1)).$
\end{thm}

\begin{thm} For $w\in S^{(N)}$, 
$[\Omega_w] = \fS_w^{(N)}(\sigma)$  in $H^*(\pflag)$.
\end{thm}
 
The quantum cohomology ring is defined as in section \ref{qmult},
with deformation variables
$q_1,\ldots,q_m$, where $q_i$ is of degree $n_{i+1}-n_{i-1}$.
For $\bd=(d_1,\ldots,d_l)$,  the corresponding 
hyperquot scheme $\R=\R(\pflag)$  
parametrizes flat families of successive quotients of $V_\Pone^*$ 
of rank $n-n_i$ and relative degree $d_{m+1-i}$.
There is a universal sequence of sheaves:
$$A_1 \hrarr A_2 \hrarr \cdots \hrarr A_{m} \hrarr V^*_{\Pone\times \R} 
\twoheadrightarrow B_{m} \twoheadrightarrow \cdots \twoheadrightarrow B_{1}$$
where $A_i$ is locally free of rank $n_i$, and $A_i\rarr A_{i+1}$ is an 
injection of sheaves, not bundles.
The subschemes $\overline{\Omega}_w$ of $\Pone\times\R$ are defined
by the appropriate degeneracy conditions.  We denote by $\mu_w$ the class of
$\overline{\Omega}_w(t)$, viewed as a subscheme of $\R$.  Details can be found
in \cite{cf2}.

\subsection{Quantum Schubert polynomials}

Let $g_i[j], c_k(l),$ and $\fS_w(c)$ be as in section \ref{usbdef}.
Define $\fS_w^{(N)}(c)$ to be the polynomial resulting from replacing
$c_i(j)$ with $c_i(n_k)$ for $n_k\leq j<n_{k+1}$.

Define two sets:
\begin{eqnarray*}
g_\sigma(N)&=&
\{(i,j):n_{l-1}+1=i\leq  j\leq n_l\mbox{ for some }1\leq l\leq m+1\}\\
g_q(N)&=& \{(i,j):1\leq i\leq j=n_l\mbox{ for some }1\leq l\leq m+1\}.
\end{eqnarray*}

Define $\E_k(l)(g)$ to be the result upon setting $g_i[j-i]=0$ for
all $(i,j)\not\in g_\sigma(N)\cup g_q(N)$ in $E_k(l)(g)$.
Define $\fS_w^{(N)}(g)$ to be the polynomial resulting from the substitution 
$c_k(n_l)=\E_k(l)$ into  $\fS_w^{(N)}(c)$.  Then the polynomial rings
generated by the $c_k(n_l)$ and the $g_i[j-i]$ are the same, so that
each $c_k(n_l)$ can be written in terms of the $g_i[j-i]$, and vice versa.

An application of the degeneracy locus formula in
 Remark 3.8 of \cite{fu}  gives
 $\mu_w = [\overline{\Omega}_w(t)]=\fS_w^{(N)}(C)=\fS_w^{(N)}(Q)$,
where $C_k(n_l)=c_k((A_l)_t^*)$ and $Q_i[j-i]$ are defined by 
the correspondence between $\Z[c]$ and $\Z[g]$.
As a purely algebraic fact, we have:
\begin{lem} \label{cg} Given a sequence of vector bundles on a scheme $X$
$$E_{m+1}\rarr E_m\rarr\cdots\rarr E_1$$
with $E_l$ of rank $n_l$ and $c_k(E_l)=\E_k(l)(g)$, we have
$$c_j(\ker(E_{l+1}\rarr E_{l}))=\left\{\begin{array}{ll}g_{n_{l-1}+1}[j-1] &\mbox{for }1\leq j<n_{l}-n_{l-1}\\
g_{n_{l+1}-j+1}[j-1] &\mbox{for }n_l-n_{l-1}\leq j\leq n_{l+1}-n_{l-1}\end{array}
\right.$$
\end{lem}

\subsection{The results} As in the case of complete flag manifolds, the
presentation of the quantum cohomology ring and quantum Giambelli formula
follow from a Main Proposition, as stated in section \ref{results},
and Lemmas \ref{polynomial} and \ref{id}, whose proofs carry through unchanged.
The difference lies in the inductive step in the
proof of Proposition \ref{q1mult}.

For the remainder of this section,
$P(g)$ denotes a polynomial in variables $g_i[j-i]$, 
for $(i,j)\in g_\sigma(N)\cup g_q(N)$, with $P^k(g)$ as in 
section \ref{mainprop}. 
Let $\D=\sum_{j=1}^m d_j(n_j-n_{j-1})$ so that $\dim \R = \dim F +\D$.
Denote
by $P(\sigma,q)$ the polynomial that results from the substitutions
\begin{align*}
g_{n_{l-1}+1}[j-1]&=\sigma^l_{j}\mbox{ when } 1\leq j\leq n_l-n_{l-1},\\
g_{n_{l-1}+1}[n_{l+1}-n_{l-1}-1]&=(-1)^{n_{l+1}-n_l-1}q_l,
\end{align*}
for $1\leq l\leq m+1$, and all other $g_i[j]=0$.  Denote
$\e^q_k(l)=\E_k(l)(\sigma,q)$.  With this notation, we have 
$\fS_w^{(N)}(\sigma,q)$ equal to the quantum
Giambelli polynomials as defined in \cite{cf2}. The results of the paper
are as follows:

\begin{thm}
$$QH^*(\pflag,\Z)\eqq \Z[\sigma_i^j,q_1,\ldots,q_m]/
(\e^q_1(m+1),\ldots,\e^q_n(m+1))$$
\end{thm}   

\begin{thm} For $w\in S^{(N)}$, 
$\sigma_w = \fS^{(N)}_w(\sigma,q)$ in $QH^*(\pflag)$.
\end{thm}

\subsection{Proof}
We first note that the proofs of the partial flag versions of
Lemmas \ref{muk+1} and \ref{Qisigma} carry through unchanged,
except for the last equality of the proof of
Lemma \ref{muk+1}, which now reads
$$\fS_w^k(\sigma,q)+\sum (-1)^{n_{i+1}-n_i+1}q_{i} 
= \fS_w^{k+1}(\sigma,q)=\fS_w(\sigma,q),$$
where the sum is over $(i,k+1)=(n_{j-1}+1,n_{j+1}-n_{j-1}-1)$.

The Main Proposition follows immediately from Proposition \ref{Pk},
whose inductive step is given by 
Lemmas \ref{muk+1} and \ref{Qisigma} and the following two propositions.
The base case is given by the same arguments as in Proposition \ref{x}.

\begin{prop} 
\label{leftpf}
Assume that Theorem \ref{Pk} holds for $k=n_{l+1}-n_{l-1}-1$.  Then
\label{qk}$$(Q_{n_{l-1}+1}[n_{l+1}-n_{l-1}+1]\cdot \mu_{w^0})_\bel 
=(-1)^{n_{l+1}-n_l+1}.$$
\end{prop}

This is the analog of Proposition \ref{q1mult} in section \ref{mainprop}. 
The proof is found in section \ref{pfleftpf}.

\begin{prop} For  $(i,j)\in g_q(N),1\leq l\leq m+1$,
$$\qsum (Q_i[j-i]\cdot\mu_w)_\bd \sw = \left\{ \begin{array}{ll}
(-1)^{n_{l+1}-n_l+1}q_l & \mbox{if } (i,j) =(n_{l-1}+1,n_{l+1}) \\
0 & \mbox{otherwise.} \end{array}\right.$$
\end{prop}
\bpf As in the complete flag case, the proof is based on dimension counts
and an understanding of bundle maps.  
Consider $(i,j)$ and $\bd$ such that the intersection number
 $Q_i[j-i]\cdot\mu_w$ in $\R$ is nonzero. By dimension considerations, 
this implies that $j-i+1\geq \D$.  Let $l,l'$ be the unique
integers so that $j=n_{l+1}$ and $n_{l'}+1\leq i<n_{l'+1}$, which implies that 
$l'\leq l+1$.  Then
$d_i\neq 0$ for $l'+1\leq i <l$, so that 
$\D=\sum d_i(n_{i+1}-n_{i-1})\geq n_{l+1}+n_{l}-n_{l'}-n_{l'-1}$.  

Since $D\leq j-i+1\leq n_{l+1}-n_{l'}$, $n_{l}-n_{l'-1}\leq 0$, 
each inequality must be an equality,
so that $l'=l-1,i=n_{l-1}+1$, and $d_l\neq 0$.  
By dimension considerations,
$\bd=\bel$ and $w=w^0$, the permutation of longest length.
Since $(w^0)^c=id$,  Proposition \ref{qk}  concludes the proof.
\epf

\subsection{Proof of Proposition \ref{leftpf}}
\label{pfleftpf}

The proof uses some of the constructions and results of \cite{cf2}, including:
\begin{equation}
h_\be^{-1}(\overline{\Omega}_w(t))= \pi^{-1}(\Pone\times \overline{\Omega}_w(t)) \cup \widetilde{\Omega}_w(t)
\end{equation}
with $\widetilde{{\Omega}}_w(t)$ of codimension $l(\tilde{w}^\be)$ in 
$\U_\be(t)$, with $l(w)-l(\tilde{w}^\be)\leq \sum e_i(n_i-n_{i-1})$.

Denote by $\A_j$ the locally free sheaf $A_j\otimes \oh_{p\times\Rl}$
for any fixed point $p\in\Pone$.
For the permutations $\alpha_{i,j}=s_{n_j-i+1}\cdots s_{n_j}$
and $\beta_{i,j}=s_{n_j+i-1}\cdots s_{n_j}$ in $S^{(N)}$, we have 
$\mu_{\alpha_{i,j}}=c_i(\A_j^*)$ and $\mu_{\beta_{i,j}}=c_i(-\A_j)$, and:
$$\tilde{\alpha}^\bel_{i,j}= \left\{\begin{array}{ll}
\alpha_{i,j} &\mbox{if } j\neq l,\\ 
\alpha_{i,l}s_{n_l}  &\mbox{if } j=l.
\end{array}\right.$$
$$\tilde{\beta}^\bel_{i,j}=\left\{\begin{array}{ll}
\beta_{i,j} &\mbox{if } j\neq l, \\ 
id &\mbox{if } j=l, 1\leq i\leq n_l-n_{l-1}, \\
\beta_{i,l}s_{n_l}\cdots s_{n_{l+1}-1} &\mbox{if } j=l,   n_l-n_{l-1}<i,
\end{array}\right.$$
where  $\bel$ is the multiindex with all zeros except a $1$
at the $i$th position. 
We have the following two lemmas:

\begin{lem} \label{chigh} For $j\geq 1$, $t\neq t'\in\Pone$,
$$h_\bel(\U_\bel)\cap \overline{\Omega}_{\alpha_{n_l-n_{l-1}+j,l}}(t) \cap \overline{\Omega}_{w^0}(t') =\emptyset.$$
\end{lem}

\begin{lem} \label{clow} For $j\geq 1$, $t_1\in\Pone$, $t_2\neq t'\in\Pone$,
$$h_\bel(\U_\bel)\cap \overline{\Omega}_{\alpha_{j,l+1}}(t_1)\cap
 \overline{\Omega}_{\alpha_{n_l-n_{l-1},l}}(t_2) 
\cap \overline{\Omega}_{w^0}(t')
=\emptyset.$$
\end{lem}

\bpf
The lemmas follow from the fact that $\dim \U_\bel = \dim F +n_l-n_{l-1}$
and dimension counts based on the following:
For any $w\in S$, $\pi^{-1}(\Pone\times \overline{\Omega}_{w}(t))$ is 
codimension $l(w)$ in $\U_\bel$ and $\widetilde{\Omega}_{w}(t)$  
is codimension $l( \tilde{w})$ in $\U_\bel(t)$.  In particular,
$\widetilde{\Omega}_{w^0}(t)$  is codimension
$\dim F-1$ in $\U_\bel(t)$ and hence codimension $\dim F$ in $\U_\bel$.
Furthermore, $\widetilde{{\Omega}}_{\alpha_{i,l}(t)}$ is codimension
 $l(\tilde{\alpha}_{i,l})=i-1$  in $\U_\bel(t')$.
Setting $i=n_l-n_{l-1}+j\geq n_l-n_{l-1}+1$ gives the proof
of Lemma \ref{chigh}.

Setting $i=n_l-n_{l-1}$,  the fact that
$\widetilde{{\Omega}}_{\alpha_{j,l+1}(t)}$ is codimension
 $l(\tilde{\alpha}_{j,l+1})=j$  in $\U_\bel(t')$ gives the proof
of Lemma \ref{clow}.
\epf

As a consequence of Lemma \ref{cg}, we have
\begin{multline}
\label{bigsum}
c_{n_{l+1}-n_{l-1}}(\A_{l+1}^*)=\sum_{j=0}^{n_l-n_{l-1}} 
Q_{n_l+1-j}[n_{l+1}-n_l-1+j]c_{n_l-n_{l-1}-j}(\A_l^*) \\
 +\sum_{j=1}^{n_{l+1}-n_{l}}
Q_{n_{l}+j+1}[n_{l+1}-n_l-j-1]c_{n_l-n_{l-1}+j}(\A_l^*).
\end{multline}
Recall that $c_k(\A_{l}^*)=\mu_{\alpha_{k,l}}$.  Then intersection of the left
hand side with $\mu_{w^0}$ in $H^*(\Rl)$ is 
zero by the argument in Lemma \ref{id}.
 For $1\leq j\leq n_l-n_{l-1}-1$, the intersection of the terms
in the first sum with $\mu_{w^0}$ are zero by the assumption on
$k$.  The intersection of the terms in the final sum with $\mu_{w^0}$ 
are zero by  Lemma \ref{chigh}. 

Since
$Q_{n_l+1}[n_{l+1}-n_l-1]=(-1)^{n_{l+1}-n_l}c_{n_{l+1}-n_l}(\A_{l+1}/\A_l)$
by Lemma \ref{cg}, we have the equality
\begin{equation}\label{expandQ}
Q_{n_l+1}[n_{l+1}-n_l-1] = \sum_{i=0}^{n_{l+1}-n_l}
(-1)^i \mu_{\alpha_{n_{l+1}-n_l-i,l+1}}
\mu_{\beta_{i,l}}.
\end{equation}
For $t\neq t'$, by Lemma 6.2 and Proposition 6.3 of \cite{cf2}, 
the intersection
$$\oO_{\alpha_{n_{l+1}-n_l-i,l+1}}(t)\cap\oO_{\beta_{i,l}}(t)
\cap\oO_{\alpha_{n_l-n_{l-1},l}}(t)\cap\oO_{w^0}(t')$$
lies in $h_\bei(\U_\bei)$ and is equal to the intersection number
$$\mu_{\alpha_{n_{l+1}-n_l-i,l+1}}\cdot\mu_{\beta_{i,l}}\cdot\mu_{\alpha_{n_l-n_{l-1},l}}\cdot\mu_{w^0}$$
in $H^*(\Rl)$.
By Lemma \ref{clow}, this is zero  for $1\leq i<n_{l+1}-n_l$.

Therefore, after substituting (\ref{expandQ}) into the intersection of
(\ref{bigsum}) with $\mu_{w^0}$, 
we get the equality
$$ 0=(Q_{n_{l-1}+1}[n_{l+1}-n_l-1]\cdot\mu_{w^0})_\bel +
(-1)^{n_{l+1}-n_l}
(\mu_{\beta_{n_{l+1}-n_l,l}}\cdot\mu_{\alpha_{n_l-n_{l-1},l}}\cdot\mu_{w^0})_\bel$$
where the last intersection can also be computed as the
Gromov-Witten number 
\begin{equation}\label{intnumber}
\langle\Omega_{\beta_{n_{l+1}-n_l,l}},\Omega_{\alpha_{n_l-n_{l-1},l}},
\Omega_{w^0}\rangle_\bel.
\end{equation}

This is a direct computation on the space of lines on $\pflag$: 
\begin{eqnarray*}
\Omega_{\beta_{n_{l+1}-n_l,l}}(A_\bullet)&=&\{L_\bullet: \dim(L_l\cap
A)\geq n_{l+1}-n_l \} \\
\Omega_{\alpha_{n_l-n_{l-1},l}}(B_\bullet)&=&\{L_\bullet: \dim( L_l\cap B)
\geq 1\}, \mbox{ and } \\
\Omega_{w^0}(C_\bullet) &=& \{L_\bullet:L_i=C_i \}\end{eqnarray*}
where $\dim L_i =\dim C_i= n_i$, $\dim A=n_{l-1}+1$, and $\dim B=n_{l+1}-1$. 

A line on $\pflag$ is of the form 
$$\{L_\bullet: L_i=D_i \mbox{ for }i\neq l, E\subset L_l \subset E'\}$$
where $D_i, E$, and $E'$ are fixed subspaces of dimension $i$, $n_l-1$, and
$n_l+1$. 
There is exactly one such line passing through these three Schubert varieties,
given by $D_i=C_i$ for $i\neq l$, $E=C_l\cap B$, and 
$E'=\langle A,C_{l}\rangle$.
Therefore, the number in the intersection (\ref{intnumber}) is one, and
$$(Q_{n_{l-1}+1}[n_{l+1}-n_l-1]\cdot\mu_{w^0})_\bel = -(-1)^{n_{l+1}-n_l}
=(-1)^{n_{l+1}-n_l+1}$$
as needed. 

This is equivalent to the single computation giving 
the relations of 
the quantum cohomology ring of the Grassmannian \cite{st}\cite{fp}.
Alternatively, by Lemma 6.2 and Proposition 6.3 of \cite{cf2}, the intersection
can be computed as the length of the zero-dimensional
scheme
$$\overline{\Omega}_{\beta_{n_{l+1}-n_l,l}}(t)\cap
 \overline{\Omega}_{\alpha_{n_l-n_{l-1},l}}(t) \cap \overline{\Omega}_{w^0}(t')
\subset h_\bel(\U_\bel),$$ or as 
$\deg([\tilde{\Omega}_{\beta_{n_{l+1}-n_l,l}}(t)]\cdot
\pi^*[t\times\overline{\Omega}_{w^0}(t')])$ in $\U_\bel$.  
We can use the projection formula, and the construction of $\U_\bel$ via a
projective bundle over $\pflag$ to obtain the result.

This concludes the proof of Proposition \ref{leftpf}, and hence
of the Main Proposition and the results of the paper.

\vspace{.5cm}
\noindent 
\sc{Department of Mathematics, Columbia University, New York, NY 10027}

\vspace{.3cm}
\noindent {\tt{lchen@math.columbia.edu}}

\end{document}